\documentclass[11pt,leqno]{amsart}
\usepackage{amsmath,amsfonts,amsthm,amscd,amssymb,graphicx}
\usepackage{epstopdf}
\numberwithin{equation}{section}
\usepackage{fullpage}
\usepackage{epstopdf}

\newcommand\cond{{($\mathcal{C}$) }}

\def\l{\lambda}
\def\eps{\varepsilon }

\newcommand\R{\mathbb R}
\newcommand\C{\mathbb C}

\def\eps{\varepsilon}

\def\l{\lambda}

\newcommand\br{\begin{remark}}
\newcommand\er{\end{remark}}
\newcommand\bp{\begin{pmatrix}}
\newcommand\ep{\end{pmatrix}}
\newcommand{\be}{\begin{equation}}
\newcommand{\ee}{\end{equation}}
\newcommand\ba{\begin{equation}\begin{aligned}}
\newcommand\ea{\end{aligned}\end{equation}}

\newcommand{\bap}{\begin{app}}
\newcommand{\eap}{\end{app}}
\newcommand{\begs}{\begin{exams}}
\newcommand{\eegs}{\end{exams}}
\newcommand{\beg}{\begin{example}}
\newcommand{\eeg}{\end{exaplem}}
\newcommand{\bpr}{\begin{proposition}}
\newcommand{\epr}{\end{proposition}}
\newcommand{\bt}{\begin{theorem}}
\newcommand{\et}{\end{theorem}}
\newcommand{\bc}{\begin{corollary}}
\newcommand{\ec}{\end{corollary}}
\newcommand{\bl}{\begin{lemma}}
\newcommand{\el}{\end{lemma}}
\newcommand{\bd}{\begin{definition}}
\newcommand{\ed}{\end{definition}}
\newcommand{\brs}{\begin{remarks}}
\newcommand{\ers}{\end{remarks}}

\newcommand{\RR}{{\mathbb R}}

\newcommand{\const}{\text{\rm constant}}
\newcommand{\Id}{{\rm Id }}

\newtheorem{theorem}{Theorem}[section]
\newtheorem{proposition}[theorem]{Proposition}
\newtheorem{corollary}[theorem]{Corollary}
\newtheorem{lemma}[theorem]{Lemma}

\theoremstyle{remark}
\newtheorem{remark}[theorem]{Remark}
\theoremstyle{definition}
\newtheorem{definition}[theorem]{Definition}

\newtheorem{example}[theorem]{Example}

\newcommand{\beq}{\begin{equation}}
\newcommand{\eeq}{\end{equation}}

\title{
Turing patterns in parabolic systems of conservation laws
and numerically observed stability of periodic waves
}

\author{Blake Barker}
\address{Brigham Young University, Provo, UT 84602}
\email{blake@math.byu.edu}
\thanks{Research of B.B. was partially supported under NSF grant no. DMS-1400872. }
\author{Soyeun Jung}
\address{Kongju National University, Korea}
\email{soyjung@kongju.ac.kr}
\thanks{Research of S.J. was partially supported by the National Research Foundation of Korea(NRF) grant funded by the
Korea government(MSIP) (No. 2016009978).
}
\author{Kevin Zumbrun}
\address{Indiana University, Bloomington, IN 47405}
\email{kzumbrun@indiana.edu} 
\thanks{Research of K.Z. was partially supported
under NSF grants no. DMS-0300487 and DMS-0801745.}
\begin{document}

\begin{abstract}
Turing patterns on unbounded domains have been widely studied in systems of reaction-diffusion equations. 
However, up to now, they have not been studied for systems of conservation laws.
Here, we (i) derive conditions for Turing instability in conservation laws and (ii) use these conditions
to find families of periodic solutions bifurcating from uniform states, numerically continuing these 
families into the large-amplitude regime.  
For the examples studied, numerical stability analysis suggests 
that stable periodic waves can emerge either from supercritical Turing bifurcations
or, via secondary bifurcation as amplitude is increased, from sub-critical Turing bifurcations.
This answers in the affirmative a question of Oh-Zumbrun whether stable periodic solutions of conservation laws can occur.
Determination of a full small-amplitude stability diagram--
specifically, determination of rigorous Eckhaus-type stability conditions--
remains an interesting open problem.
\end{abstract}
\date{\today}
\maketitle


\section{Introduction}

The study of periodic solutions of conservation laws and their stability, initiated in \cite{Oh2003a,Oh2003} and
continued in \cite{Se,JZ2}, etc., has led to a number of interesting developments, particularly in the related study of 
roll-waves in inclined shallow-water flow.  
For an account of these developments, see, e.g., \cite{Johnson2012a} and references therein.
However, in the original context of conservation laws, so far {\it no example of a stable periodic wave has been found.}
Indeed, one of the primary results of \cite{Oh2003a,Pogan2013} was that for the fundamental example of planar viscoelasticity,
stable periodic waves do not exist, due to a special variational structure of this particular system;
it was cited as a basic open problem whether stable periodic waves could arise for {\it any} system of conservation laws,
either physically motivated: or artificially contrived.

In the more standard context of reaction diffusion systems and classical pattern formation theory, by contrast,
stable periodic solutions are abundant and well-understood, through the mechanism of {\it Turing instability}, or
bifurcation of small-amplitude, approximately-constant period, periodic solutions from a uniform state.
For such waves, stability is completely determined by an associated {\it Eckhaus stability diagram},
as derived formally in \cite{Eckhaus1965} and verified rigorously in \cite{Mielke, Mielke1997, Sc2,SZJV},
essentially by perturbation from constant-coefficient linearized behavior.
By contrast, the small-amplitude waves investigated up to now (see Remark \ref{btrmk})
come through more complicated zero-wave number bifurcations in which period goes to infinity as amplitude goes to zero
and the stability analysis is far from constant-coefficient
(see, e.g., \cite{Barker2014} in the successfully-analyzed case of shallow-water flow).

Our simple goal in this paper, therefore, is {\it to seek stable periodic waves via a conservation law analog of Turing
instability.}
In the first part, we find an analog of Turing instability, 
with which we are able to generate large numbers of examples of spatially periodic solutions of conservation laws.
Next, we find an interesting dimensional restriction to systems of three or more coordinates, explaining
the absence of Turing instabilities for $2\times 2$ systems considered previously.
Finally, we perform a numerical existence/stability study for $3\times 3$ example systems exhibiting Turing instability,
answering in the affirmative the fundamental question posed in \cite{Oh2003a,Pogan2013} whether there can exist
stable spatially periodic solutions of systems of conservation laws, at least at
the level of numerical approximation. 
These studies suggest that, for supercritical Turing bifurcation, stable waves can emerge through the small-amplitude
limit and persist up to rather large amplitudes.  For sub-critical Turing bifurcations, all emerging waves
are necessarily initially unstable, but appear in some cases to undergo secondary bifurcation to stability as
amplitude is further increased.

The numerically observed stability of intermediate-amplitude waves we regard as conclusive.
Delicacy of numerical approximation as amplitude goes to zero, however,
prevents us from obtaining a detailed stability diagram near the Turing bifurcation or 
even from making definitive conclusions about stability in that regime.
Rigorous spectral stability analysis for conservation laws in this regime, 
analogous to those of \cite{Mielke, Mielke1997, Sc2,SZJV}
in the reaction diffusion case, we regard therefore as a very interesting open problem.
The studies in \cite{Matthews2000,Sukhtayev2016} of reaction diffusion equations with
an associated conservation law may offer guidance in such an investigation.


\section{Turing instability for conservation laws}\label{s:turing}

We begin by defining a notion of Turing instability for systems of conservation laws
\be
\label{sys}
u_t + f(u;\eps)_x=(D^\eps u_x)_x,
\ee
$u\in \RR^n$, where $\eps$ is a bifurcation parameter and $D$ for simplicity is taken constant.
Linearizing \eqref{sys} about a uniform state $u(x,t)\equiv u_0$ yields the family of constant-coefficient equations
\be\label{lin}
u_t=L(\eps)u:=  - A^\eps u_x + D^\eps u_{xx}
\ee
with dispersion relations $\lambda_j(\xi)\in \sigma(-i\xi A^\eps -\xi^2D^\eps)$, $\xi\in \R$, where
$\sigma(\cdot)$ here and elsewhere denotes spectrum of a matrix or linear operator.
The state $u_0$ is spectrally (hence nonlinearly) stable if
\be\label{con1}
\Re \sigma(-i\xi A^\eps-\xi^2D^\eps )\leq -\theta |\xi|^2, \quad \theta>0,
\ee
for all $\xi\in \R$ \cite{Kaw}.

Following the original philosophy applied by Turing \cite{T} to reaction diffusion systems, we seek a natural
set of conditions guaranteeing {\it low- and high-frequency stability}-- i.e., 
that \eqref{con1} hold for $|\xi|\to 0, \infty$-- but allowing instability at finite frequencies
$|\xi|\neq 0, \infty$.
Should this be possible, then performing a homotopy in $\eps$ between stable and unstable states, we
may expect generically to arrive at a special bifurcation point $\eps=\eps_*$, without loss of generality $\eps_*=0$,
for which \eqref{con1} holds uniformly away from special points $\xi=\pm \xi_* $, at which
\be\label{max}
\max_{\xi\neq 0} \Re \sigma(-i\xi A^{\eps_*}-\xi^2D^{\eps_*})=0
\ee
is achieved (note, by complex conjugate symmetry, that extrema appear in $\pm$ pairs) and for which \eqref{con1} fails strictly
as $\eps$ is further increased.
We may then conclude, by standard bifurcation theory applied to the domain of periodic functions with period 
$X:=2\pi/\xi_*$ the appearance of nontrivial spatially periodic solutions with periods near $X$,
similarly as in the reaction diffusion case \cite{Mielke, Mielke1997, Sc2,SZJV}.

At $\xi=0$, \eqref{con1} yields that $A$ is hyperbolic, in the sense that it has real semisimple eigenvalues.
Without loss of generality, therefore, take $A$ diagonal, with entries $a_j$, $j=1, \dots, n$.
In the simplest case that $A$ is {\it strictly hyperbolic,} in the sense that these $a_j$ are distinct, we find by
spectral perturbation expansion about $\xi=0$ \cite{Kaw} that the corresponding eigenvalue expansions are
$$ 
\lambda_j(\xi)=-ia_j 
\xi-D_{jj}\xi^2 + O(\xi^3),
$$
so that \eqref{con1} ($\xi\ll 1$) is equivalent to the condition that $D$ {\it have positive diagonal
entries $D_{jj}$}.
Similarly, by spectral expansion about $\xi=\infty$, 
$$
\sigma(-i\xi A-\xi^2D)=-\xi^2 \sigma(D) + O(\xi),
$$
so that \eqref{con1}($\xi=\infty$) is equivalent to the condition that $D$ be unstable, i.e.,
{\it have eigenvalues with strictly positive real part}.
Collecting, our conditions are \cond:

$\bullet$ $A^\eps$ is diagonal with distinct entries, and 

$\bullet$ $D^\eps$ has positive diagonal entries and eigenvalues with strictly positive real part.
\medskip

\noindent
These are to be contrasted with Turing's conditions in the reaction diffusion case $u_t=Du_{xx} +g(u)$
that $D$ be symmetric positive and $A:=dg(u)$ be symmetric negative definite \cite{T}.

\subsection{Turing instability and Hopf bifurcation}
Let \eqref{max} hold, with 
$\lambda= \pm i\tau \in \sigma(-i\xi A-\xi^2 D)$ for $ \xi=\pm \xi_*$,
$\xi_*\ne 0$.
Then, changing to the moving coordinate frame $ x\to  \tilde x:=x-ct$, for 
$ c:=\tau/\xi_*, $ 
or, equivalently, under the the change 
of coordinates $A\to \tilde A:=A-cI$, we have
$\lambda= 0 \in \sigma(-i\xi \tilde A-\xi^2 D)$ for $ \xi=\pm \xi_*$,
i.e., 
$ \det( -i\xi \tilde A-\xi^2D)=0
$
at $\xi=\xi_*$, or
\be\label{Hopf}
\pm i\xi_*\in \sigma (D^{-1}\tilde A).
\ee

Condition \eqref{Hopf} may be recognized as the condition for Hopf
bifurcation of an equilibrium $u(x,t)\equiv \const$ of the traveling-wave ODE
\be\label{ode}
D^\eps u'=f(u;\eps)-cu + q,
\ee
where $q$ is a constant of integration, for which the linearized equation is
$ u'=D^{-1}\tilde A u, $ $\tilde A$ again diagonal.
Thus, we recover by finite-dimensional bifurcation theory the previously-remarked appearance of nontrivial
periodic solutions with period near $X=2\pi/\xi_*$.
We also obtain the alternative bifurcation criterion \eqref{Hopf}.
{\it This simplifies the problem a great deal;} for one thing, we are now
working with real matrices, as occur for
symbols in the reaction diffusion case, and not complex ones.

\subsubsection{Dimensional count}\label{s:dim}
From the usual Hopf bifurcation theorem for ODE, we find that for each fixed nearby $q$, $c$, 
there exists a one-parameter family of nontrivial periodic solutions bifurcating from the constant solution, generically
parametrized nonsingularly by period $X$.
Thus, fixing $q=0$, we obtain a 
{\it $2$-parameter family of periodic solutions, generically well-parametrized by $c$ and $X$.}

\subsection{Finding Turing instabilities}
To find Turing instability, we may seek $A^\eps$ and $D^\eps$ satisfying \cond, $\eps\in \R$ a bifurcation parameter,
such that \eqref{max} is violated at $\eps=0$ (instability), but \eqref{con1} is satisfied for all $\xi$ at $\eps=1$ (stability),
for example if $A^1=\Id$ or $D^1=\Id$.
For, in this case, the conditions \cond on $A^\eps$, $D^\eps$ insure that at the largest value 
$\eps_*$ of $\eps$ for which \eqref{max} is satisfied, the maximum \eqref{max} is achieved at some $\xi=\xi_*\neq 0$,
while for $\eps>0$ there must be strictly positive real part eigenvalues, again bounded uniformly away from zero.

As another approach, starting from the observation relating Turing instabilities and Hopf bifurcation, 
notice first that \eqref{max} cannot occur for $D=I$, in which case the spectra
of $(-i\xi A-\xi^2D)$ are simply $\lambda_j(\xi)=-i\xi a_j - \xi^2$; nor can \eqref{Hopf}, since $\sigma(\tilde A)$ is by
assumption real.
Thus, we suggest, first, finding examples $\check A$, $\check D$ satisfying
\eqref{Hopf} either analytically or by checking random matrices,
then, setting up a homotopy $D^\epsilon:= \epsilon \check D + (1-\epsilon)I$
from the identity to $\check D$.
Since, as just observed, $\sigma(-i\xi \check A-\xi^2 D^\epsilon)$ is stable
for $\eps=0$, while for $\eps=1$ it is at most neutrally stable, having
zero eigenvalues at $\xi=\pm \xi_*\neq 0$, we find that for some
$\epsilon \in (0,1]$,
$\sigma(-i\xi \check A-\xi^2 D^\epsilon)$ is exactly neutral, i.e., 
a {\it Turing instability}, with eigenvalues $\pm i\tau$ at $\xi=\pm \hat \xi_*$ (note: different from the original $\xi_*$ in general!).
As described above, this corresponds to a Hopf bifurcation in the traveling-wave
ODE for speed $c_*:= \tau/\hat \xi_*$, with limiting wave number $\hat \xi_*$
and period $X_*:=2\pi/\hat \xi_*$.

\section{Negative results}\label{s:negative}

We next describe situations in which Turing instability {\it cannot} occur, narrowing our search.

\subsection{The $2\times 2$ case}
We have the following result for $n=2$, strikingly different
from the situation of the
reaction diffusion case.

\bpr\label{noprop}
Assuming \cond, there exist \emph{no} Turing-type instabilities of \eqref{sys} for $n=2$.
\epr

\begin{proof}
Take by assumption $A$ diagonal.
Since $D^{-1}A$ is real, appearance of a pure imaginary eigenvalue $i\tau$
implies the appearance also of its complex conjugate $-i\tau$,
hence trace is zero and determinant is positive.  By a scaling transformation
$S=\bp \alpha & 0 \\ 0 & \beta\ep$ not affecting diagonal form of $A$, we
may arrange therefore that
$
D^{-1}A=\bp c & 1\\-1 & -c\ep=:J,
$
for some $c^2<1$.
Noting that $ J^2= (c^2-1) I,  $
we may solve to obtain
$ D= \frac{1}{c^2-1}AJ= \frac{1}{c^2-1}\bp a_1c & a_1\\ -a_2 & -a_2c\ep.  $
The requirement that $D$ have positive diagonal implies, with $c^2<1$,
that $a_1 c<0$ and $a_2 c>0$, so that $a_1$ and $a_2$ have opposite sign.
But, $\det D=(c^2-1)^{-2} a_1a_2(1-c^2)>0$ implies that
$a_1 $ and $a_2$ have the same sign, hence these two conditions cannot
hold at once.
\end{proof}

\begin{example}\label{btrmk}
The viscoelasticity model
$\tau_t-u_x=d_{11} \tau_{xx}$, $u_t + p(\tau)_x= d_{22} u_{xx}$
studied by Oh-Zumbrun \cite{Oh2003a} falls into the above framework, hence does not
admit Turing instabilities.  In fact, periodic waves arise in this
model through Bogdanov-Takens bifurcation associated with splitting
of two or more equilibria, a more complicated bifurcation far from
constant-coefficient behavior.
\end{example}
\subsection{Simultaneous symmetrizability}
Another case in which Turing instabilities do not occur is when
$A$ and $D$ are simultaneously symmetrizable, or, equivalently,
can be converted by change of coordinates to be both symmetric.
For, then, in the new coordinates, $D$,
being symmetric positive definite, has a square root, and so
$D^{-1}A$ is similar to the symmetric matrix
$D^{1/2}D^{-1}AD^{-1/2}= D^{-1/2}AD^{-1/2}$, hence has {\it real
eigenvalues.}
More generally, it is easy to see that Turing instability does not
occur for $A$ symmetric and $\Re D:=(1/2)(D+D^T)>0$, since
$D^{-1}Av=i\tau v$ would imply
$
0=\Re i\tau \langle v, Av\rangle= Re \langle v, D v\rangle= \langle v, \Re D v\rangle>0, 
$
a contradiction.
This recovers the well-known fact that existence of a viscosity-compatible
convex entropy for the system \eqref{sys} implies nonexistence of 
non-constant stationary solutions, since existence of such an entropy implies
the corresponding symmetry conditions on the linearized equations.
Thus, taking $A$ without loss of generality diagonal, we must specifically
seek $D$ {\it nonsymmetric}, $D+D^T$ {\it nonpositive} in order to find Turing instability.

\subsection{Nonstrict hyperbolicity}
Finally, we give a simple example showing that the condition of strict hyperbolicity of $A^\eps$
is necessary in \cond. 
Consider the matrices 
\be\label{singeg}
A^\eps
=\begin{pmatrix} 1 & 0 & 0\\ 0 & \eps & 0 \\ 0 & 0 & 1\end{pmatrix} \quad \text{and} \quad D=\bp 1 & 0& 2 \\ 0 & 1 & 1 \\ 1 & -2 & 1 \ep. 
\ee
Here, $\sigma (D)=\{1\}$; so $-i\xi A - \xi^2 D$ is stable for $|\xi| \rightarrow +\infty$. For $|\xi| \rightarrow 0$, we look at $2\times 2$ blocks corresponding to the 1 and 3 entries of $A$ and $D$, 
\be
\tilde A = \bp 1&0 \\0 &1 \ep \quad \text{and} \quad \tilde D = \bp 1 & 2 \\ 1 & 1 \ep.
\ee
Then, the two eigenvalues of $-i\xi A - \xi^2 D$ close to $i\xi$ for $\xi\ll 1$ 
are by standard spectral perturbation theory $\l_j(\xi)=-i\xi-\xi^2 \tilde d_j$ , 
where $\tilde d_j$ are eigenvalues of $\tilde D$. We easily see that $\tilde D$ has two real eigenvalues 
with opposite sign because $det(\tilde D)=-1<0$. Thus, \eqref{con1} is not satisfied for $|\xi| \rightarrow 0$. 

\br\label{singrmk}
Though example \eqref{singeg}, failing \cond, does not itself yield Turing instability, it is quite useful
in finding nearby systems that do. For, note perturbation in $\eps$ generates matrices $D^{-1}A$ with nonstable eigenvalues 
despite $A>0$.  Perturbing first $\eps$ to obtain instability, then $A$ still more slightly to recover strict hyperbolicity,
we thus obtain an example satisfying \cond with unstable $D^{-1}A$, which yields a Turing bifurcation upon
homotopy $D\to I$.
We in fact used this method to generate the examples of Section \ref{s:num}.
(We have generated other examples in other ways, that were not reported here; all
exhibited similar behavior, however.)
\er

\section{Spectral and nonlinear stability}\label{s:framework}
Before describing our numerical investigations, we briefly recall the abstract stability framework 
developed in \cite{Oh2003a,JZ2,Johnson2012a}, etc., relevant to stability of
the nontrivial periodic waves bifurcating from a constant solution at Turing instability.
First, recall \cite{JZ2,Johnson2012a} that, under the condition of transversality of the associated periodic orbit of the
traveling-wave ODE (guaranteed in this case by the Hopf bifurcation scenario, for sufficiently small-amplitude waves),
nonlinear stability with respect to localized perturbations of the periodic wave considered as a solution on the whole line
is determined (up to mild nondegeneracy conditions) by conditions of {\it diffusive spectral stability,}
as we now describe.

By Floquet theory, the $L^2(\R)$ spectrum of the linearized operator $L$ about a periodic wave of period $X$ is 
entirely essential spectrum, corresponding to values $\lambda \in \C$ for which there exist
generalized eigenfunction solutions $v(x)=e^{i\xi x}w(x)$, $\xi\in \R$, 
of the associated eigenvalue equation $(L-\lambda)v=0$ with $w$ periodic, period $X$.
The dissipative stability conditions are that this spectrum have real part $\leq -\eta \xi^2$, $\eta>0$,
for all $\xi\in \R$, and strictly negative for $(\xi, \lambda) \neq (0,0)$.

For transversal orbits with $\eps$ bounded away from $\eps_*$, 
the spectra near $(\xi,\lambda)=(0,0)$ consists of the union of $(n+1)$ smooth spectral curves 
$ \lambda_j(\xi)=-ia_j\xi +o(\xi) $ 
through the origin $\lambda=0$, which, under the nondegeneracy condition that $a_j$ be distinct, are analytic in $\xi$,
admitting second-order expansions
\be\label{specexp}
\lambda_j(\xi)=-ia_j \xi-b_j\xi^2 +O(\xi^3), \quad j=1,\dots, n+1.
\ee
Moreover, the functions $\lambda_j(\xi)$ correspond to the linearized dispersion relations for the associated
second-order {\it Whitham system},
an associated second-order $(n+1)\times(n+1)$ system of conservation laws formally governing slow modulational 
behavior \cite{Wh,Se,Johnson2012a}.  Thus, low-frequency diffusive spectral stability is equivalent to well-posedness
(hyperbolic-parabolicity) of the Whitham system, which is in turn equivalent to reality of $a_j$ (hyperbolicity)
and positivity of $\Re b_j$ (parabolicity) in \eqref{specexp}, with high-frequency spectral stability
given by $\Re \lambda \leq -\eta<0$ for $|\xi|\geq \eta$, $\eta>0$.

In the case of Turing instability, choosing the period $X_*$ such that the wave-numbers $\pm\xi_*$ at $\eps=\eps_*$
are equal to zero modulo $2\pi/X_*$, we find by direct Fourier transform calculation that the constant solution at $\eps=\eps_*$
has low-frequency spectrum consisting of {\it $(n+2)$ spectral curves passing through the origin,}
with all other spectra satisfying $\Re \lambda \leq -\eta<0$ for some $\eta>0$.
The spectra of the bifurcating periodic waves perturbs smoothly from these values as $\eps$ is increased, hence high-frequency
diffusive stability is guaranteed.
However, low-frequency stability is now determined by a possibly complicated bifurcation of $(n+2)$ spectral curves involving the
$(n+1)$ ``Whitham curves'' \eqref{specexp} passing through the origin plus an additional curve originating from the constant limit passing close to but not through the origin.
These curves are clearly visible in the numerically approximated spectra displayed below in Section \ref{s:num} for example systems
with $n=3$: namely, $4$ Whitham curves passing through the origin, with a $5$th (initially) neutral spectral curve passing near the origin, with all $5$ of these passing through the origin at the bifurcation point $\eps=\eps_*$.

\section{Numerical investigations}\label{s:num}
Guided by the results of Sections \ref{s:turing}, \ref{s:negative}, and \ref{s:framework},
we now perform the main work of the paper,
carrying out numerical existence and stability investigations for 
periodic solutions of systems of conservation laws arising through Turing bifurcation from the uniform state 
in dimension $n=3$.
Numerics are carried out using the MATLAB-based package STABLAB developed for this purpose \cite{STABLAB}.

\subsection{Quadratic nonlinearity} We first consider the system 
\be \label{sysdef2}
u_t + A^{\eps}u_x+N(u)_x = Du_{xx},
\ee
with
\be 
A^{\eps}:= \bp 1&0&0\\0&a_{22}^0+\eps&0\\0&0&3 \ep, \quad D:= \bp1&0&2\\ 0&1&1\\1&-2&1\ep,\quad \text{and} \quad N(u):= \beta \bp u_1^2\\0\\0 \ep,
\ee
where $a_{22}^0 = 2.605173614560316$. Here,  $\eps$ is a bifurcation parameter that we will vary and $u \equiv 0$ is a constant solution of \eqref{sysdef2}. By linearization of \eqref{sysdef2} about  $u=0$, we have 
\be \label{linearsysdef2}
u_t + A^{\eps}u_x = Du_{xx}. 
\ee

\begin{figure}[htbp]
	\begin{center}
		$
		\begin{array}{lcr}
		  (a) \includegraphics[scale=0.25]{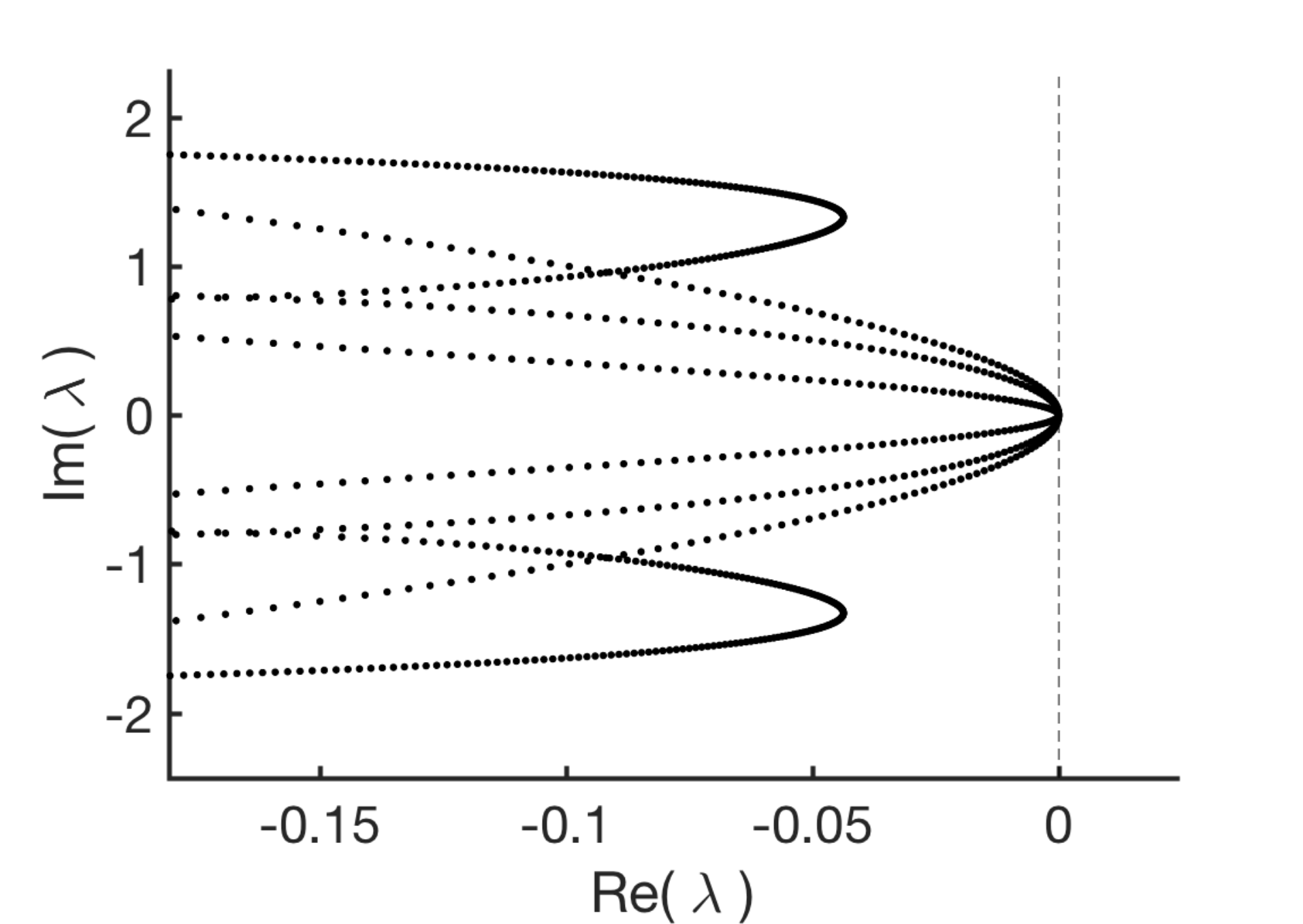}
		& (b) \includegraphics[scale=0.25]{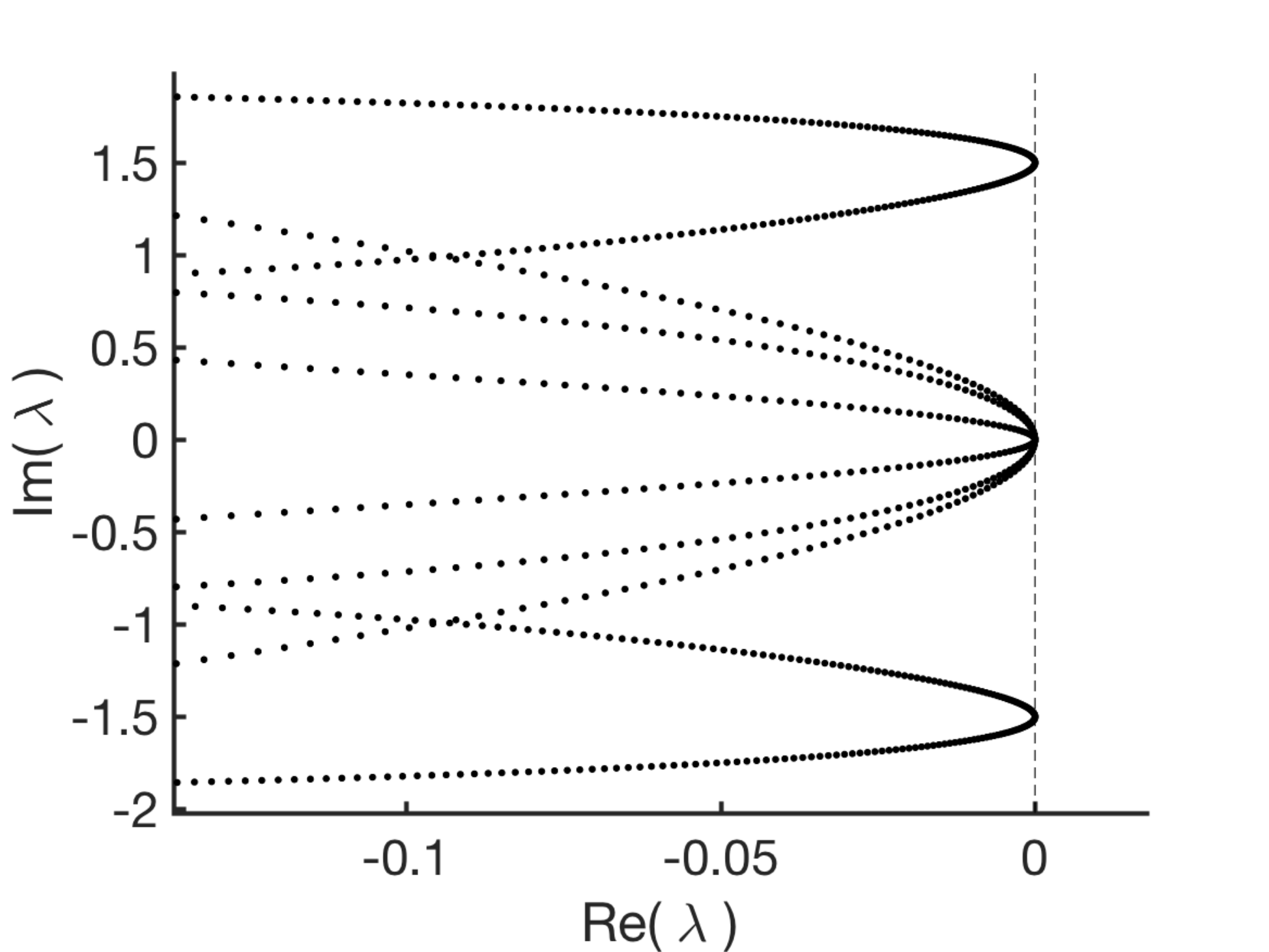}
		& (c) \includegraphics[scale=0.25]{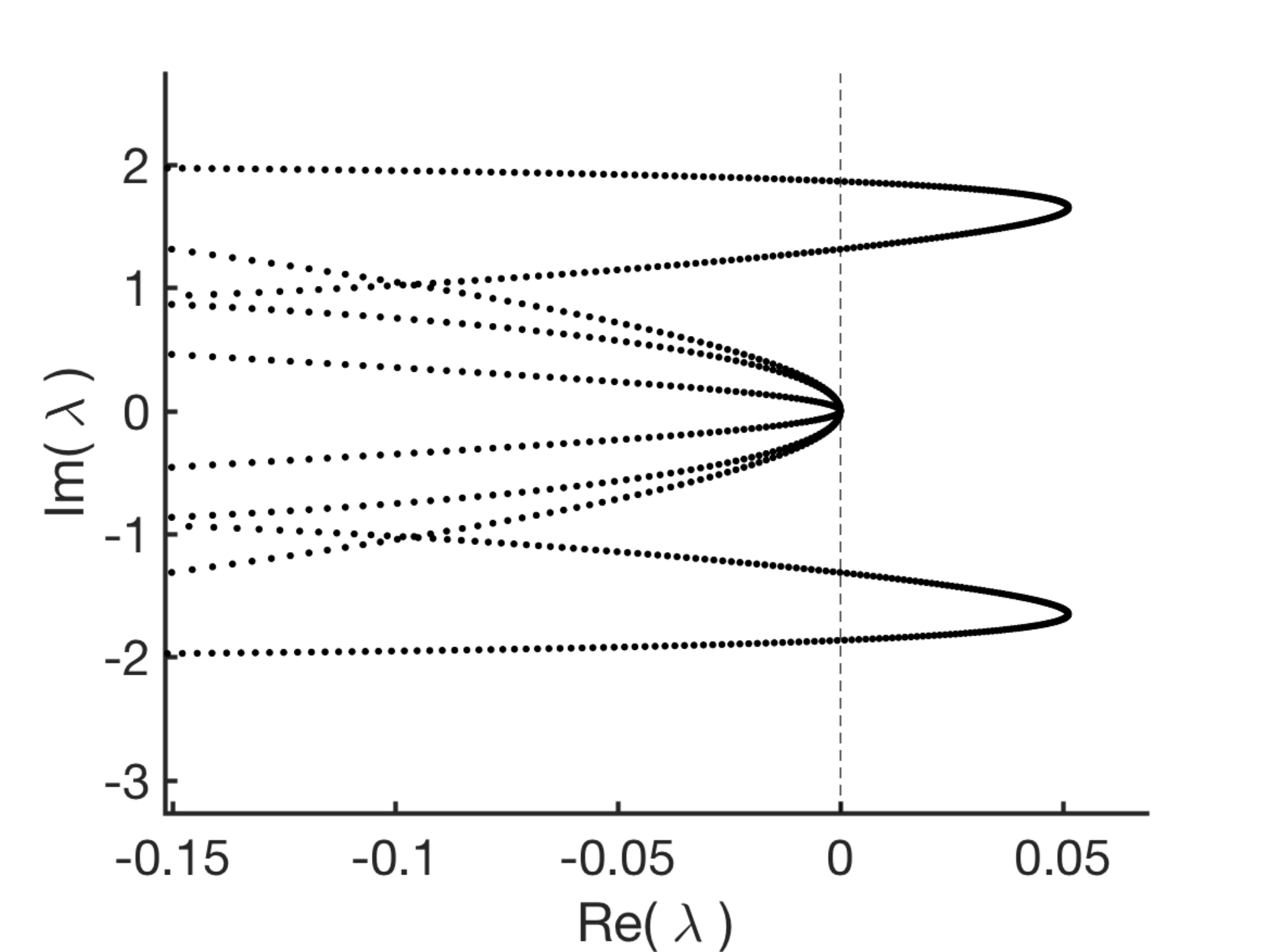}
		\end{array}
		$
	\end{center}
	\caption{Plot with dots of a sampling of the spectrum of the constant solution, $-i\xi A -\xi^2D$, with (a) $\eps = -0.2$, (b) $\eps = 0$, (c) $\eps = 0.2$. The dashed vertical line marks the imaginary axis.   
	}
	\label{super spectrum of constant solution with zero c}
\end{figure}

\begin{figure}[htbp]
	\begin{center}
		$
		\begin{array}{lcr}
		(a)\includegraphics[scale=0.25]{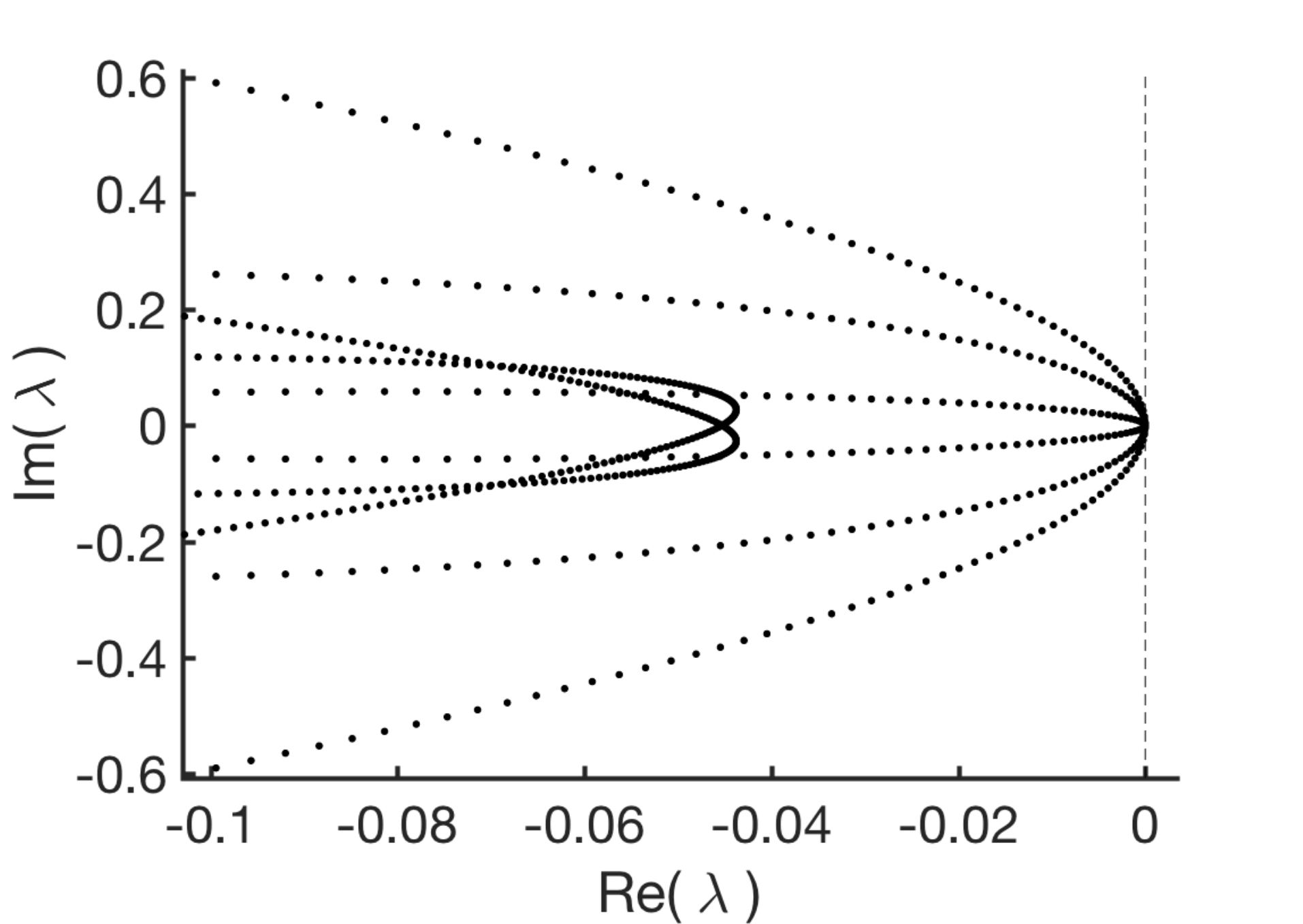}& (b)\includegraphics[scale=0.25]{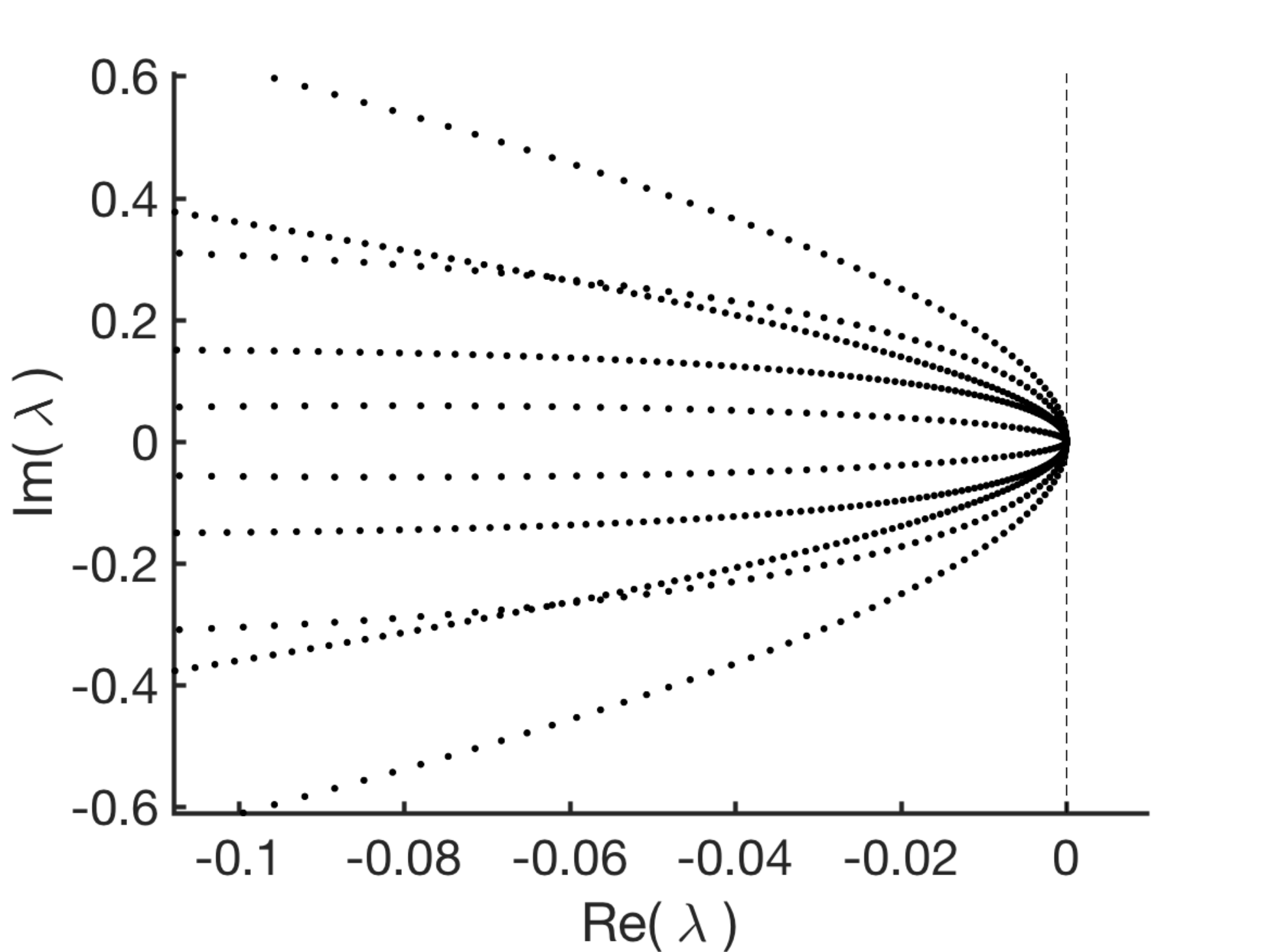}&
		(c)\includegraphics[scale=0.25]{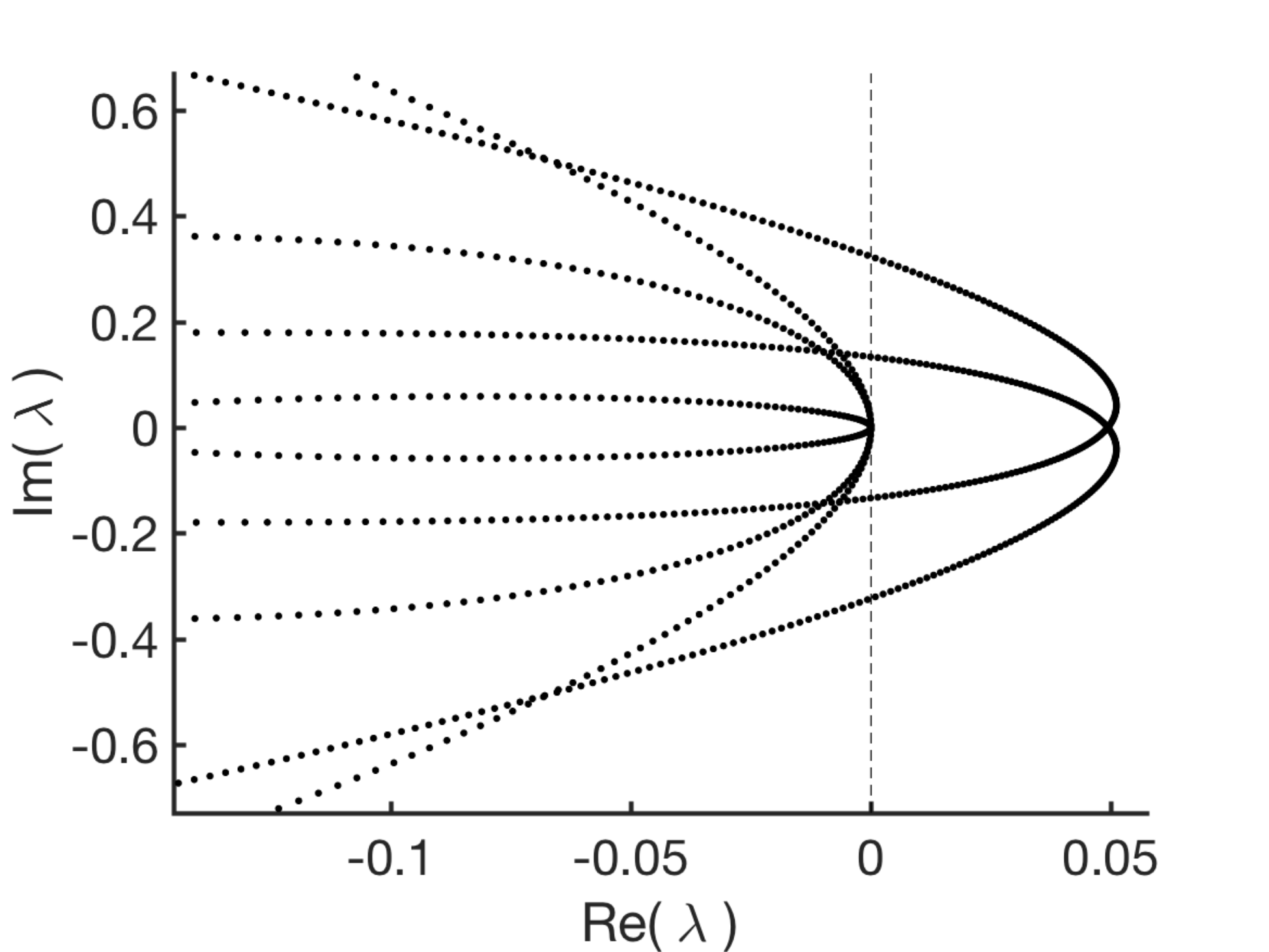}
		\end{array}
		$
	\end{center}
	\caption{Plot with dots of a sampling of the spectrum of the constant solution, $-i\xi (A-c_*I) -\xi^2D$, with (a) $\eps = -0.2$, (b) $\eps = 0$, (c) $\eps = 0.2$ and $c=c_* \approx 1.30$. The dashed vertical line marks the imaginary axis.   
	}
	\label{super spectrum of constant solution with nonzero c}
\end{figure}

We first check Turing-type instability conditions for $u \equiv 0$ in \eqref{linearsysdef2}. Notice that $A^{\eps}$ is strictly hyperbolic and $D$ has positive diagonal entries with $\sigma(D)=\{1\}$, which means that $-i\xi A^{\eps}-\xi^2 D$ is stable near
$\xi =0$ or $\xi=\pm \infty$.  We examine numerically 
stability of $u \equiv 0$ as $\eps$ changes. In Figure \ref{super spectrum of constant solution with zero c}, we plot the spectrum of $-i\xi A^{\eps}-\xi^2 D$ with $\eps = -0.2$, $\eps = 0$, and $\eps = 0.2$.  It is seen that the constant solution $u \equiv 0$ is stable for $\eps < 0$ and unstable for $\eps >0$. Thus, Turing instability occurs at $\eps =0$, that is, \eqref{max} is satisfied with $\pm i\tau \in \sigma(-i\xi A^0-\xi^2 D)$ for $\tau \approx 1.5$ and $ \xi_* \approx \pm 1.16$. As we observed in the previous section, $\pm i \xi_*$ are eigenvalues of $D^{-1}(A^0-c_* I) $ for $c_* = \frac{\tau}{\xi_*} \approx 1.30$. So the condition for Hopf bifurcation of a constant solution $u \equiv 0$ of the profile equation 
\be
-cu + A^{\eps}u+N(u) = Du'+q
\ee 
is satisfied at the bifurcating point $\eps=0$ and $c=c_*$. Here $q \in \RR^3$ is an integration constant and we fix $q=0$ from now on.  In Figure \ref{super spectrum of constant solution with nonzero c}, we plot the spectrum of $-i\xi (A^{\eps}-c_* I) - \xi^2 D$ for the same $\eps$ as in Figure \ref{super spectrum of constant solution with zero c}, showing
how this moves the neutral spectrum from $\lambda=\pm i\tau$ to $\lambda=0$.

The Hopf bifurcation leads to periodic profiles bifurcating from the uniform state $u \equiv 0$. 
In order to solve for these profiles, we let $\eps$ be a free variable and vary the period $X$ and wave speed $c$,
approximating associated solutions using the periodic profile solver built into STABLAB, which uses MATLAB's Newton-based 
boundary-value problem solver bvp5c. 
In addition to periodic boundary conditions, the profile solver specifies a phase condition $w\cdot f(y(0)) = 0$ where $y'(x) = f(y(x))$ is the profile ODE (\eqref{ode} in the present case)
and $w$ is a random vector. Unless $w$ is a degenerate choice, $w\cdot \dot y(t) = 0$ for some $t$ by periodicity of $y$ and Rolle's Theorem, so this phase condition chooses a solution (at least locally) uniquely. 
To numerically solve the profile equation with a quadratic nonlinearity, we first obtain a solution by using as an initial guess $u(x) = \sqrt{\eps}\Re(e^{2\pi i x}v )/10$, where $v$ is the real part of an eigenvector, whose corresponding eigenvalue has non-zero imaginary part, of the profile Jacobian evaluated at the fixed point $(0,0,0)^T$. 
That is, we start with an initial guess consisting of a strategically scaled periodic solution of
the linearized equations at the bifurcation point $\eps=0$.
Once we have a profile solution via this guess, we use continuation to solve for other profiles
with nearby period $X$ and speed $c$, obtaining thereby a full $2$-parameter family of approximate solutions
parametrized by $(c,X)$, as described in Section \ref{s:dim}.

\begin{figure}[htbp]
	\begin{center}
		$
		\begin{array}{lcr}
		(a) \includegraphics[scale=0.2]{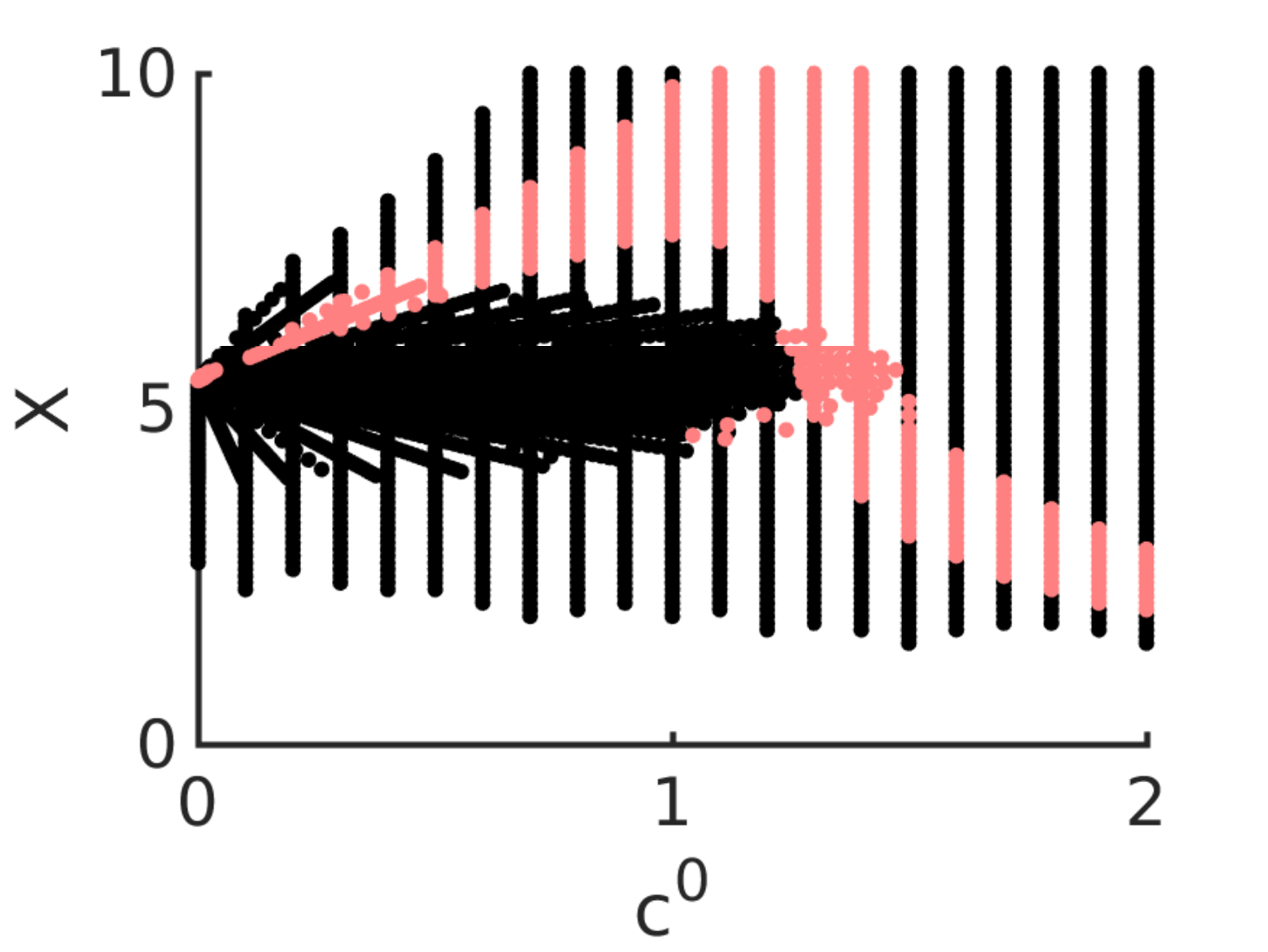}&
		(b) \includegraphics[scale=0.2]{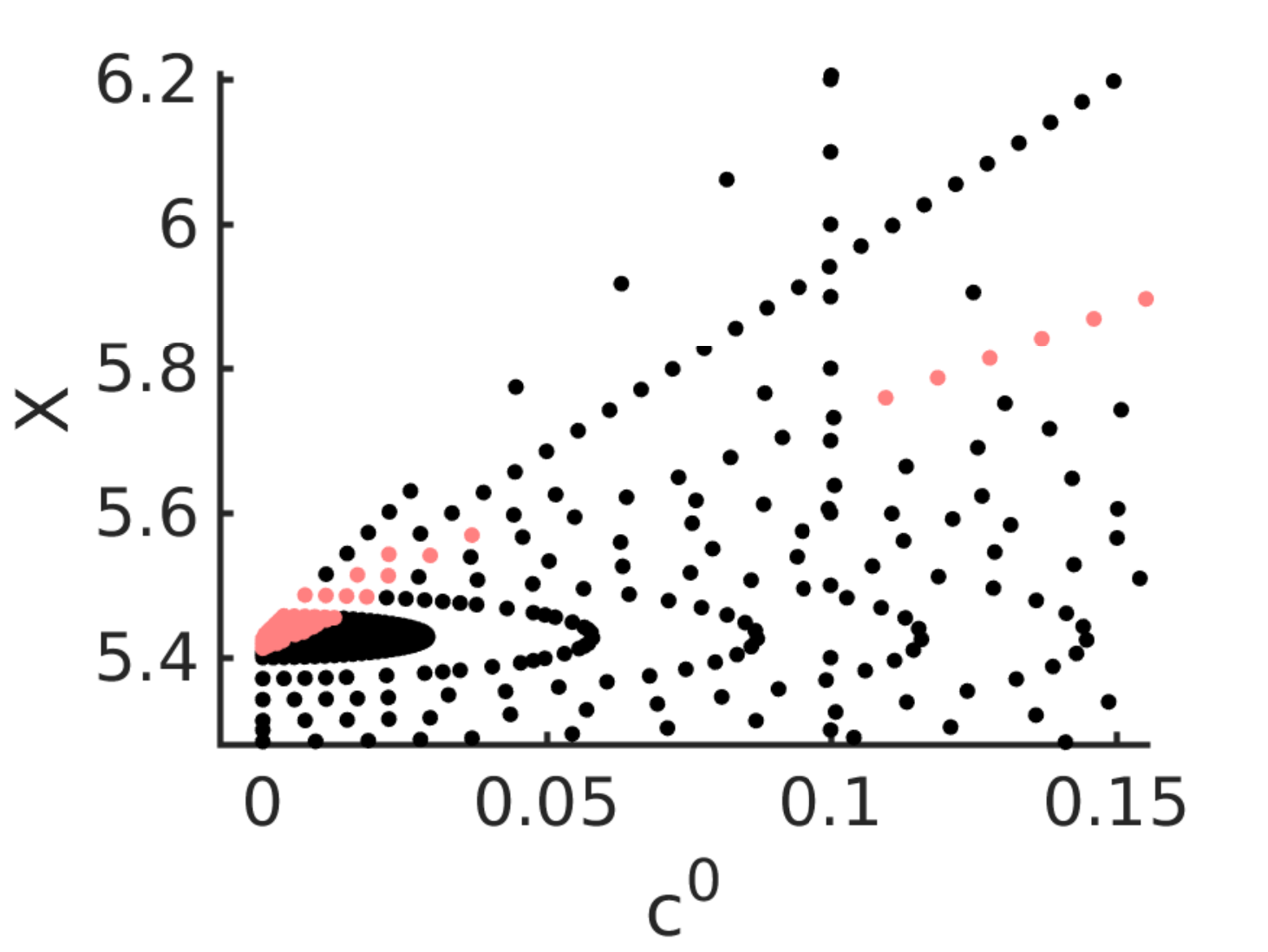} &
		(c) \includegraphics[scale=0.2]{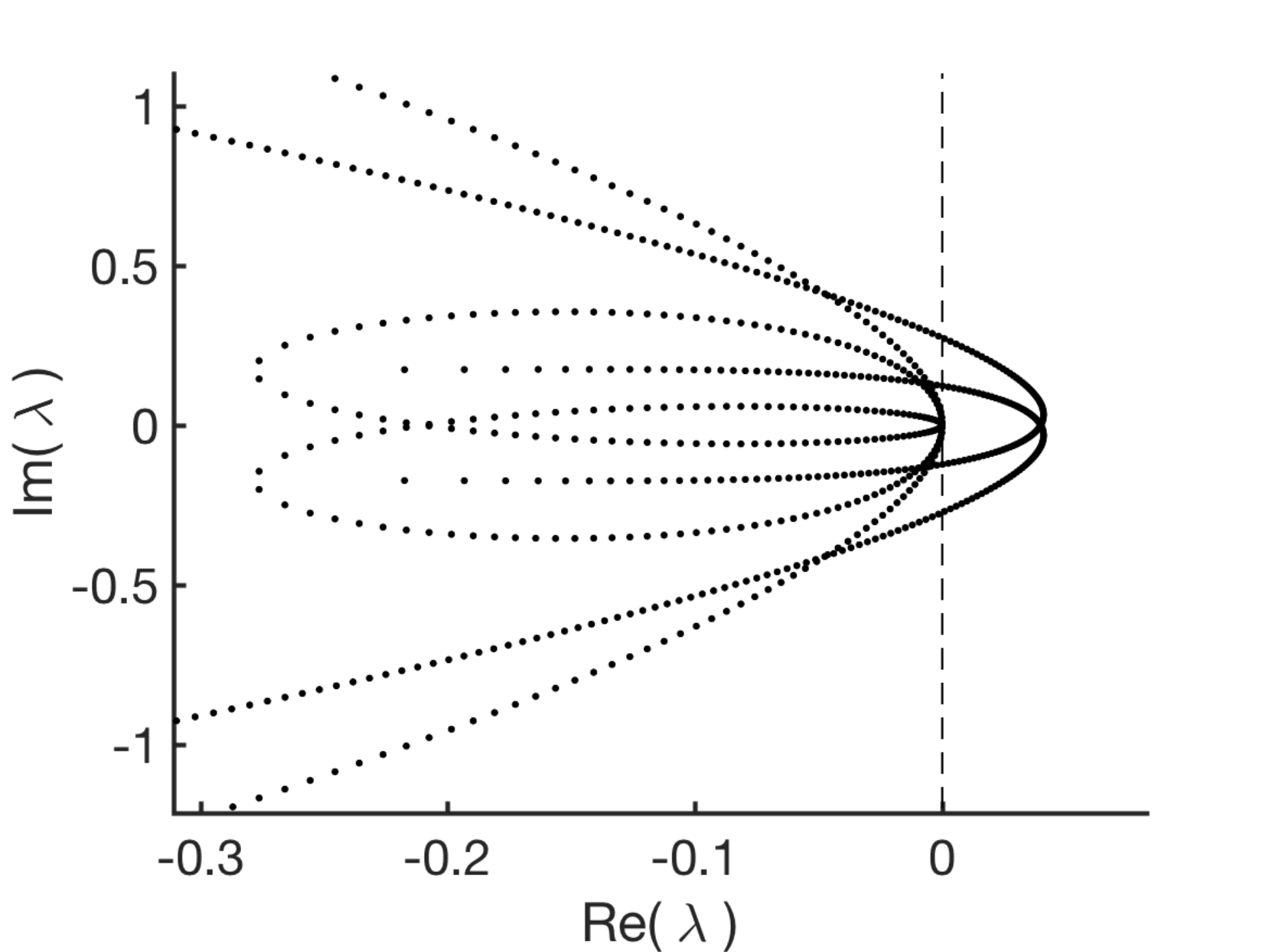}\\
		(d) \includegraphics[scale=0.2]{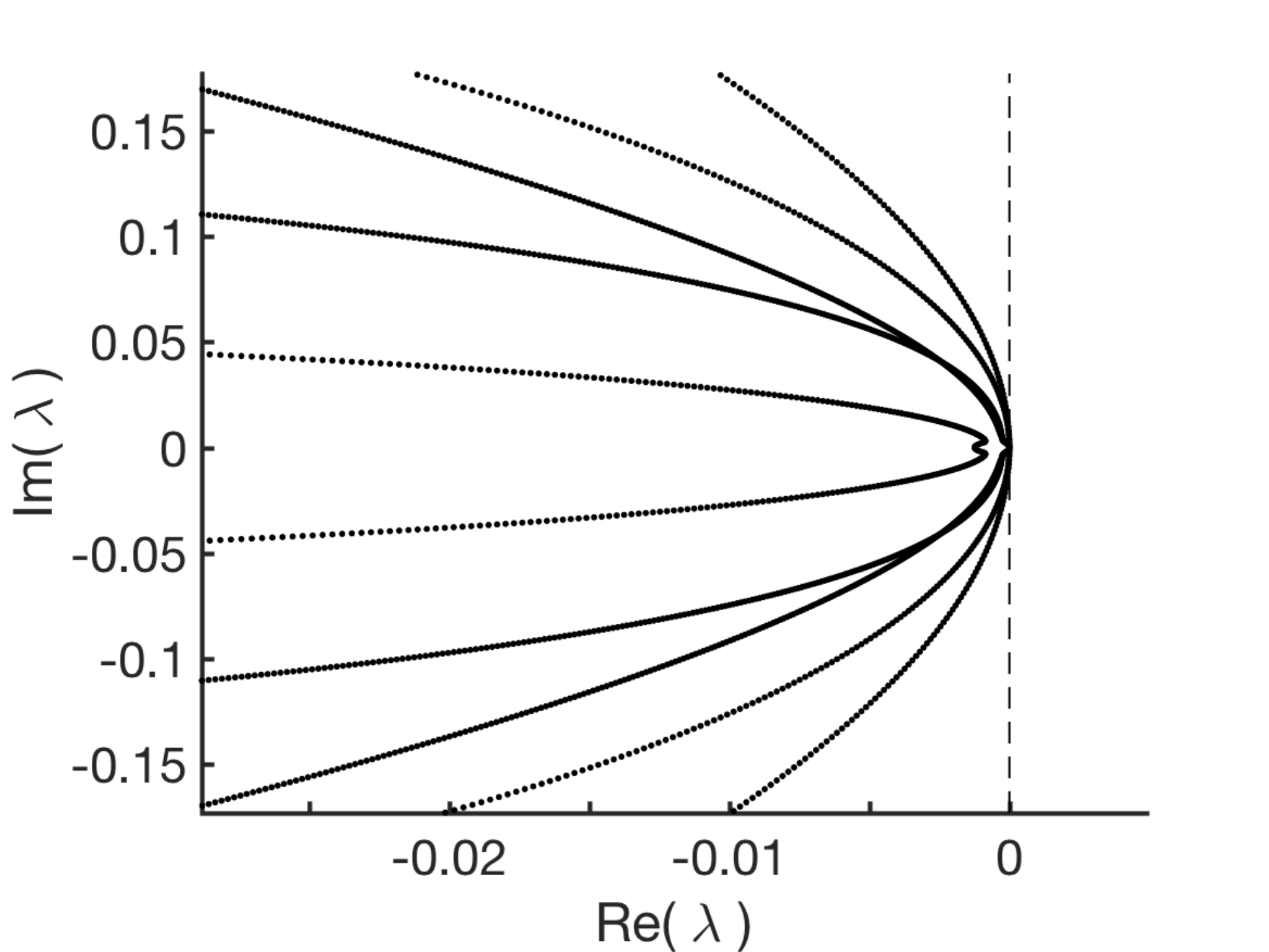}&
		(e) \includegraphics[scale=0.2]{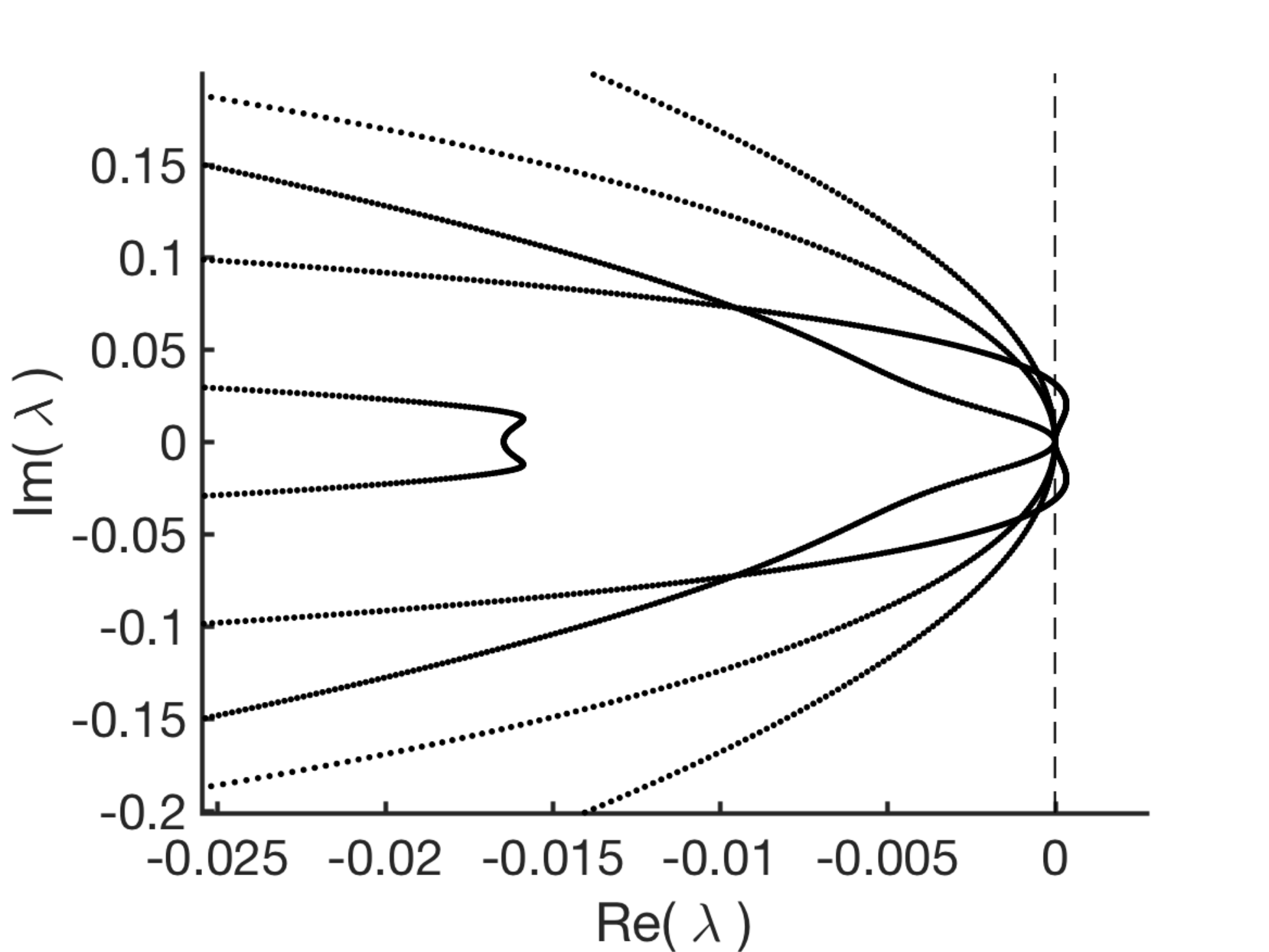}&
		(f) \includegraphics[scale=0.2]{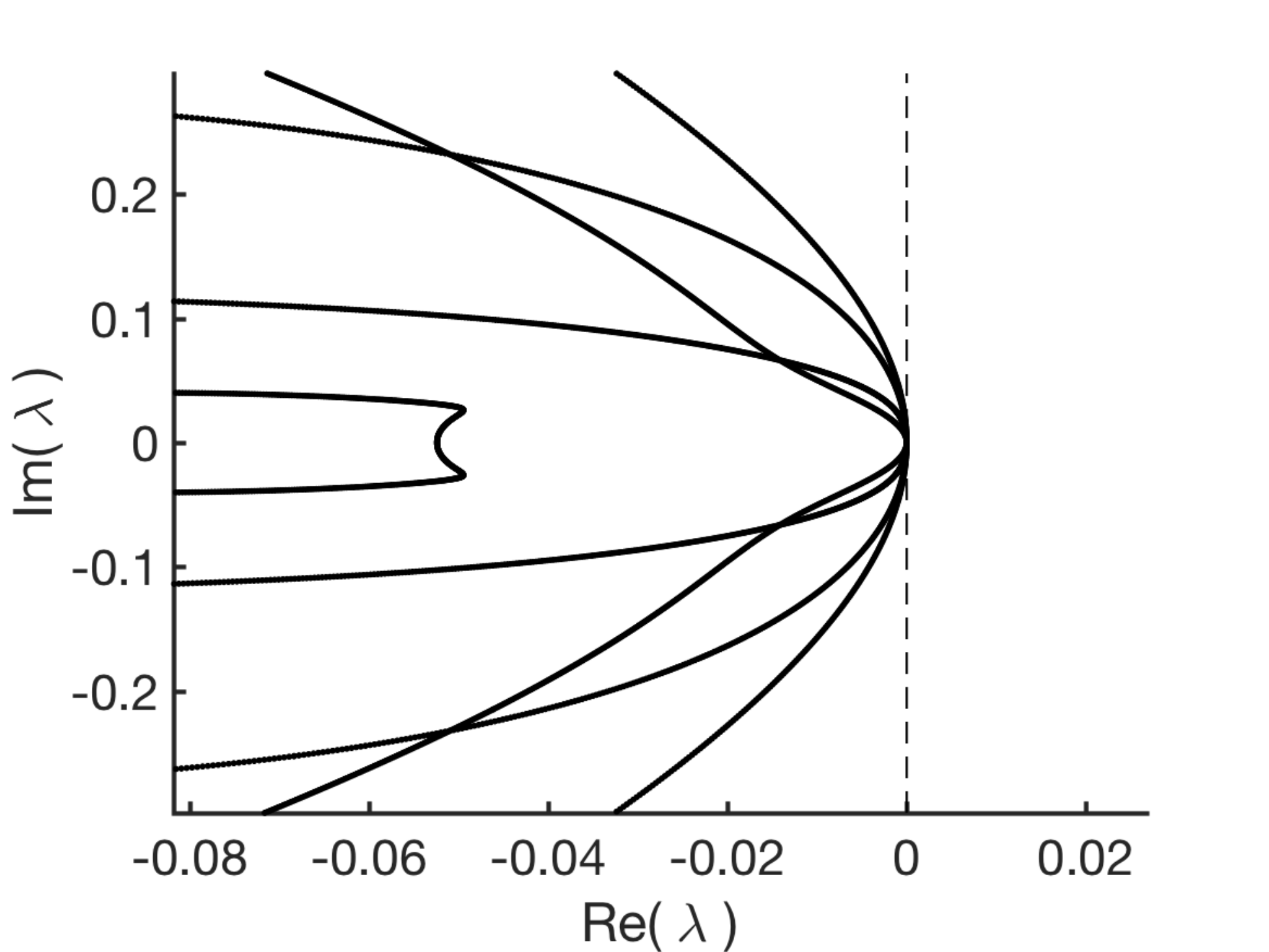}\\
		& (g) \includegraphics[scale=0.2]{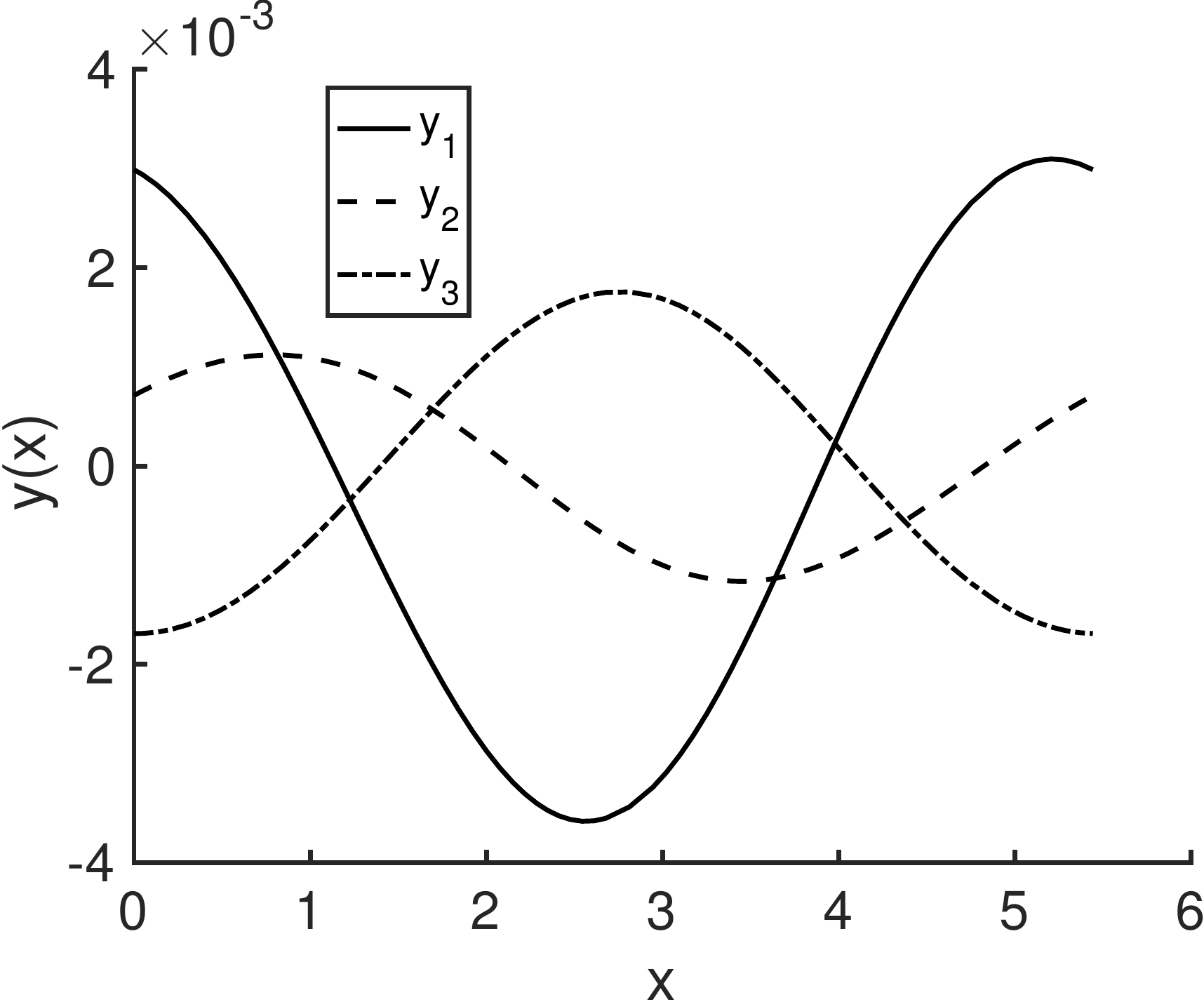} &
		\end{array}
		$
	\end{center}
	\caption{(a) Stability bifurcation diagram in the coordinates of shifted wave speed $c^0 = c-c_*$ and period $X$. Pink dots (light dots in grayscale) and black dots correspond respectively to stable and unstable waves. (b) Zoom in of (a) showing a family of stable waves in parameter space leading to the point of the Turing bifurcation. There is a small region of instability separating the stable waves near the Turing bifurcation point and the large stability region. (c) Plot of the spectrum of the zero constant solution when $\eps = 2.82e-3$, $c = c_*+4.06e-3$, and $X = 5.44$, indicating that the Turing bifurcation corresponds to a supercritical Hopf bifurcation. (d) Plot of the spectrum of a periodic wave in the family of stable waves bifurcating from the Turing bifurcation. (e) Plot of the spectrum of a periodic wave in the family of unstable waves separating the two regions of stability. (f) Plot of the spectrum of a periodic wave in the large stability region. (g) Plot of the bifurcating periodic profile when $\eps = 2.82e-3$, $c = c_*+4.06e-3$, and $X = 5.44$, with component one marked with a solid line, component two with a dashed line, and component three with a dot-dashed line.   Throughout $\beta = -10$ and a dashed line marks the imaginary axis.}
	\label{fig113}
\end{figure}

In Figure \ref{fig113} (a) and (b),  we plot the stability bifurcation diagram in the coordinates of shifted wave 
speed $c^0=c-c_*$ and period $X$. 
The bifurcation diagram shows that there is a family of stable waves bifurcating from the Turing bifurcation. There is a small region of instability occurring from a ``parabolic'' Whitham instability, or change in curvature of a neutral spectral 
curve through the origin, corresponding to negative diffusion or ill-posedness of the associated formal slow modulation 
Whitham equations, which separates the region of stability near the Turing bifurcation point and the larger stability region. Figures \ref{fig113} (d)-(f) demonstrate this onset of Whitham-type instability as seen in the spectrum of the 
bifurcating periodic waves. In Figure \ref{fig113} (c), we see that the spectrum of the background constant solution becomes unstable as $\eps$ increases, so that the periodic profile shown in Figure \ref{fig113} (g) comes into existence through a super-critical Hopf bifurcation. Finally, in Figure \ref{fig113} (g), we plot the periodic profile for $\beta = -10$, $\eps = 2.82e-3$, $c = c_*+4.06e-3$, $X = 5.44$.

We note that, as described in Section \ref{s:framework}, there are generically 4 neutral spectral curves passing through the origin, with second-order Taylor expansions related to the linearized dispersion relation for a formal Whitham slow-modulation approximation. This is clearly visible in Figure \ref{fig113} (d)-(f). However, as seen in Figure \ref{super spectrum of constant solution with nonzero c} (b), the constant solution has 5 spectral curves passing through the origin at the bifurcation point and the spectra of bifurcating periodic waves perturbs from these 5 curves. So, at the bifurcation point, there is a 5th neutral curve passing through the origin, which remains nearby for values of $\eps$ nearby $\eps_*$. It explains why the the spectrum of stable periodic waves bifurcating from Turing bifurcation in Figure \ref{fig113} (d) has an additional 5th curve which is very close to the origin but not through the origin. Stability of small-amplitude waves is determined by behavior of these 5 neutral curves, either by movement of the maximum real part of the 5th curve into the unstable or stable half-plane (``co-periodic'' stability, corresponding with super- or sub-criticality of the associated Hopf bifurcation), or by a ``Whitham-type'' instability consisting of loss of tangency to the imaginary axis (first-order, or ``hyperbolic'' instability) or change in curvature (2nd order, or ``parabolic'' instability) of one of the 4 neutral curves through the origin; see Section \ref{s:framework}.

For the quadratic nonlinearity, if $ u(x) $ is a profile solution for a fixed $ \beta $, then $ - u(x) $ 
is a profile solution for $ - \beta $, with the same value of $ \eps $. 
Thus, we are not able to produce a corresponding sub-critical Hopf bifurcation 
by reversing the sign of $\beta$, but a mirror super-critical bifurcation.

To find examples of stable periodic profiles corresponding to both sub and super-critical Hopf bifurcations, we change the quadratic nonlinearity to a cubic nonlinearity in the next example, removing this symmetry and allowing us to change from super- to sub- by changing the sign of $\beta$.


\subsection{Cubic nonlinearity} We consider next the system of conservation laws
\be \label{sysdef3}
u_t + A^{\eps}u_x+N(u)_x = Du_{xx},
\ee
with
\be 
A^{\eps}:= \bp 1&0&0\\0&a_{22}^0+\eps&0\\0&0&3 \ep, \quad D:= \bp1&0&2\\ 0&1&1\\1&-2&1\ep,\quad \text{and} \quad N(u):= \beta \bp u_1^3\\0\\0 \ep,
\ee
where $a_{22}^0 = 2.605173614560316$. Similarly as the quadratic example, we vary $\eps$ as a bifurcation parameter. 
The stability of $u \equiv 0$ as $\eps$ varies is already shown in Figure \ref{super spectrum of constant solution with zero c} and Figure \ref{super spectrum of constant solution with nonzero c}. 

Starting from the super-critical periodic profile solutions found previously
for the quadratic nonlinearity, we obtain a solution for the cubic nonlinearity by continuation in a homotopy variable $0\leq h\leq 1$ via the nonlinearity $N(U)=[\beta(hy_1^3+(1-h)y_1^2),0,0]^T$. To obtain a sub-critical profile solution for the cubic nonlinearity, we use the
approximate symmetry $(\beta,c,\eps)\to(-\beta,-c,-\eps)$, which is valid at the linear periodic level only. 
Thereafter, we solve for profiles using continuation.

\begin{figure}[htbp]
	\begin{center}
		$
		\begin{array}{lcr}
		(a) \includegraphics[scale=0.25]{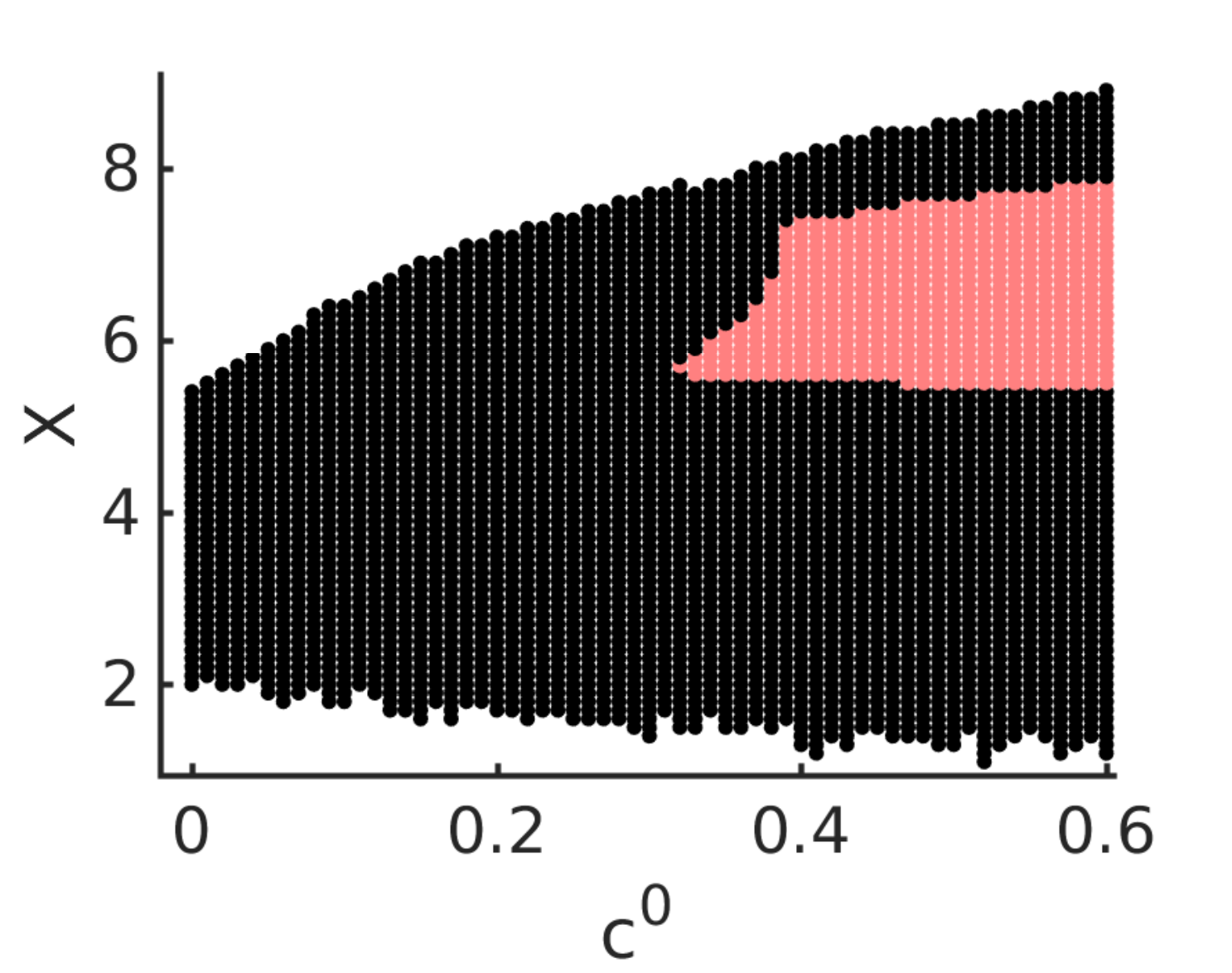} & (b) \includegraphics[scale=0.25]{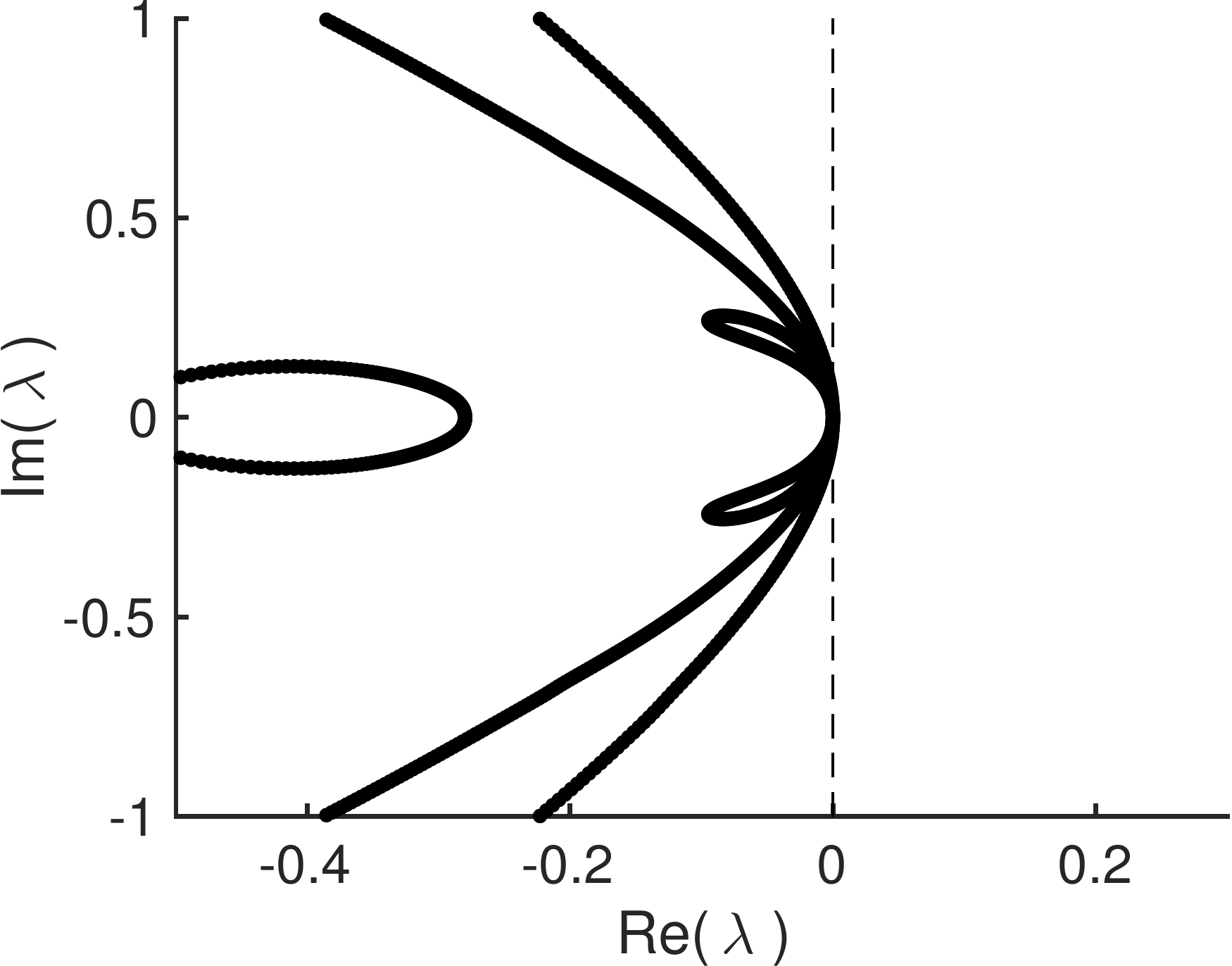} &
		(c) \includegraphics[scale=0.25]{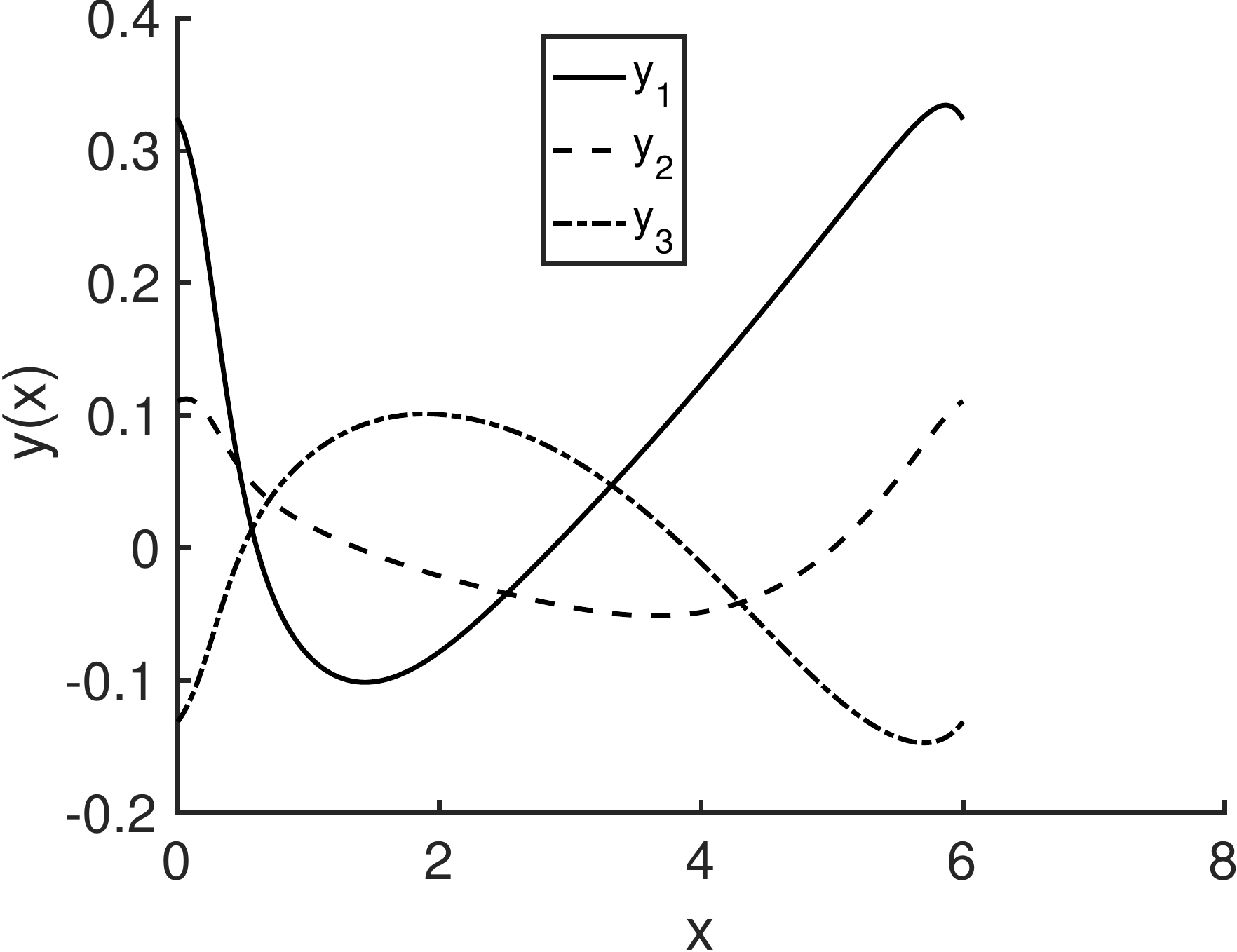}
		\end{array}
		$
	\end{center}
	\caption{(a) Stability diagram in the coordinates of shifted wave speed $c^0 = c-c_*$ and period $X$ for $\beta = 10$. Pink dots (light dots in grayscale) and black dots correspond respectively to stable and unstable waves. (b) For a stable wave, we plot in (b) its spectrum and in (c) the wave itself, with $\beta = 10$, $c^0 = 0.5$, $X = 6$ , and $\eps = 8.74e-1$. A dashed line marks the imaginary axis in (b).
	}
	\label{fig133}
\end{figure}

In Figure \ref{fig133}, we plot the bifurcating stable periodic solution through a super-critical Hopf bifurcation. Since $\eps > 0$ for the constant solution to be unstable, as seen in Figure \ref{super spectrum of constant solution with nonzero c}, the periodic profile shown in Figure \ref{fig133} (c) exists through a super-critical Hopf bifurcation. Figure \ref{fig133} (b) shows the stable spectrum of the periodic profile shown in (c). Here $\beta = 10$, $c^0 = 0.5$, $X = 6$ , and $\eps = 8.74e-1$. In Figure \ref{fig133} (a), we plot a stability diagram in the coordinates of shifted wave speed $c^0 = c-c_*$ and period $X$. We do not find a family of stable waves bifurcating from the Turing instability.

\begin{figure}[htbp]
	\begin{center}
		$
		\begin{array}{lcr}
		(a) \includegraphics[scale=0.3]{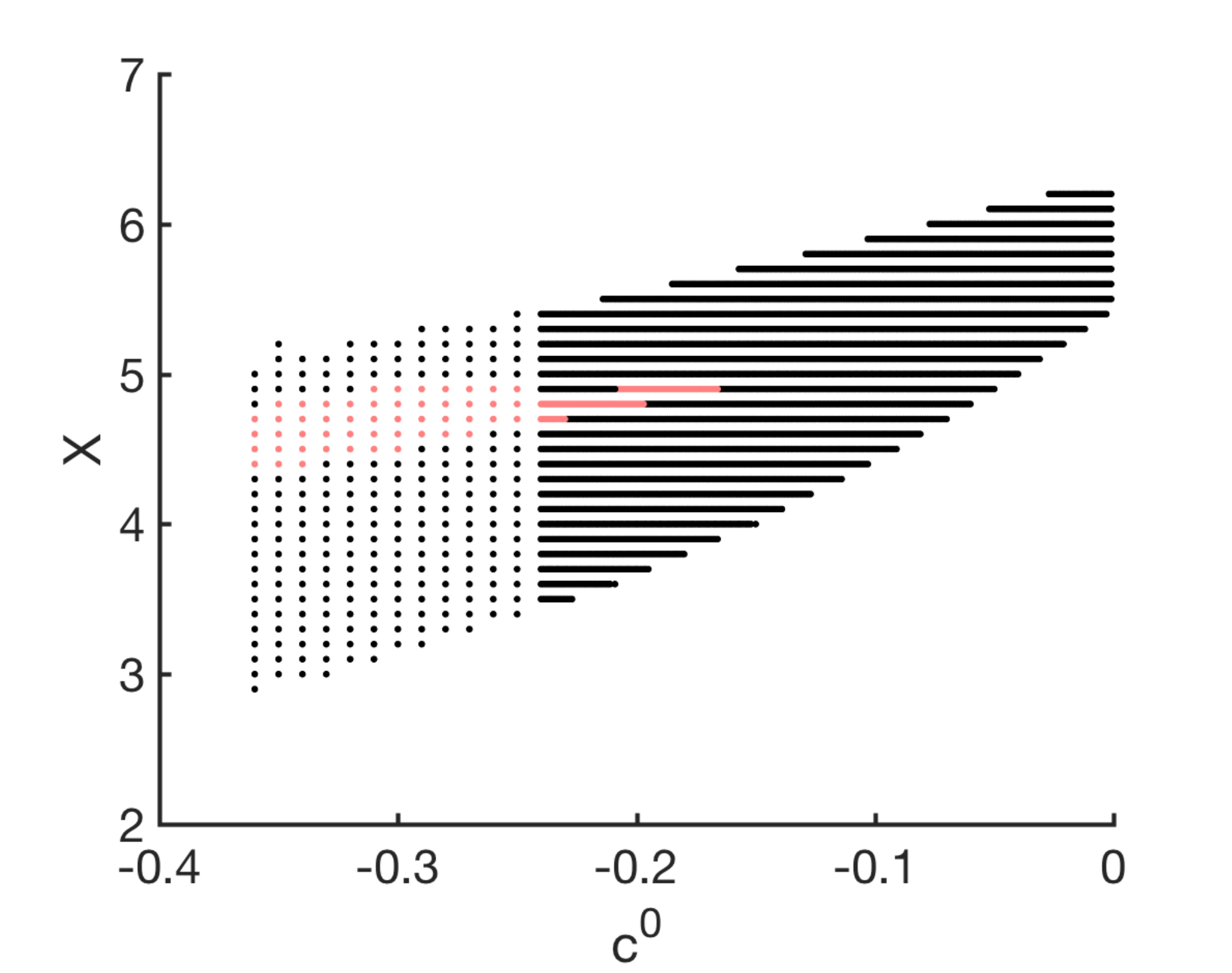} & (b) \includegraphics[scale=0.3]{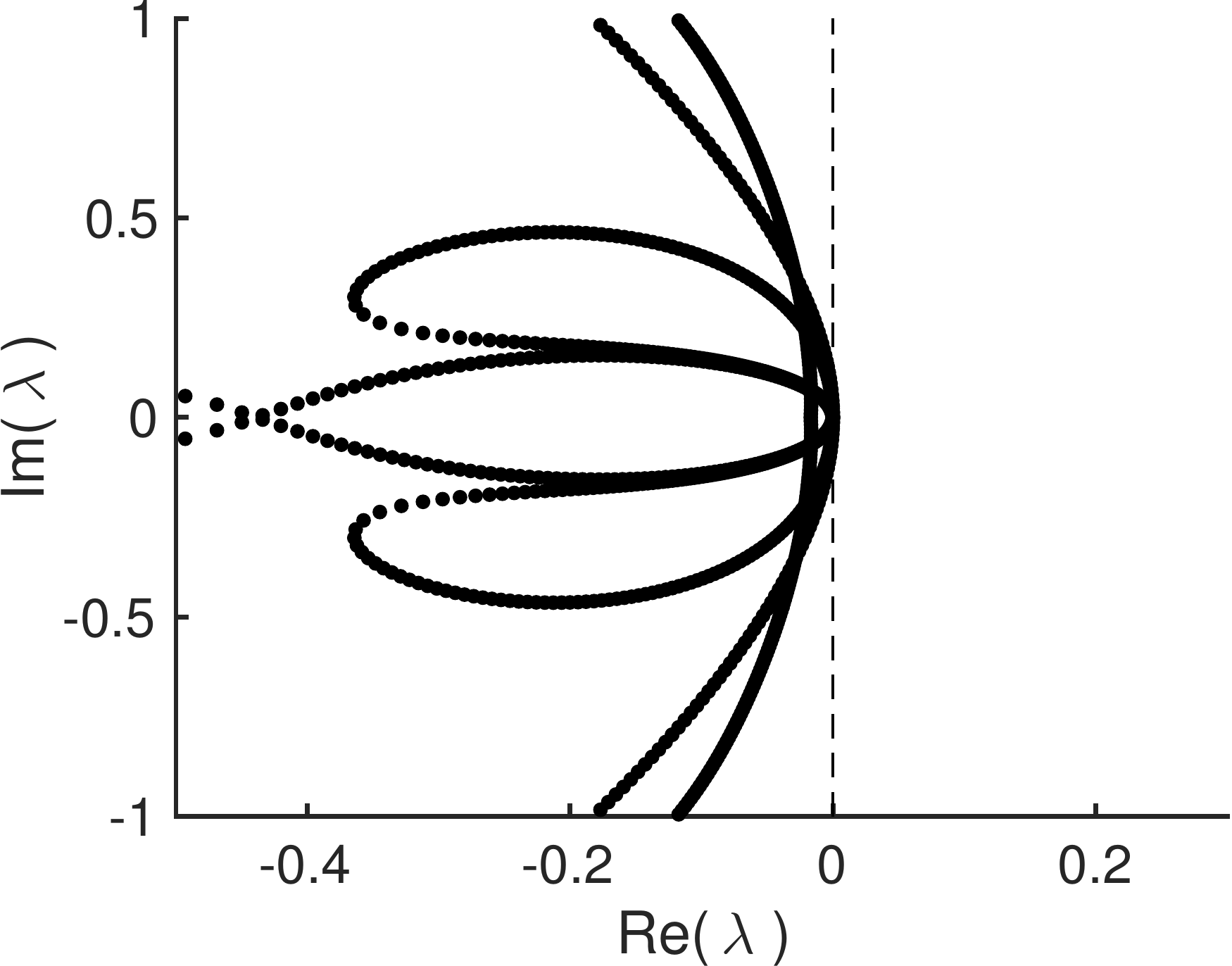}\\
		 (c) \includegraphics[scale=0.3]{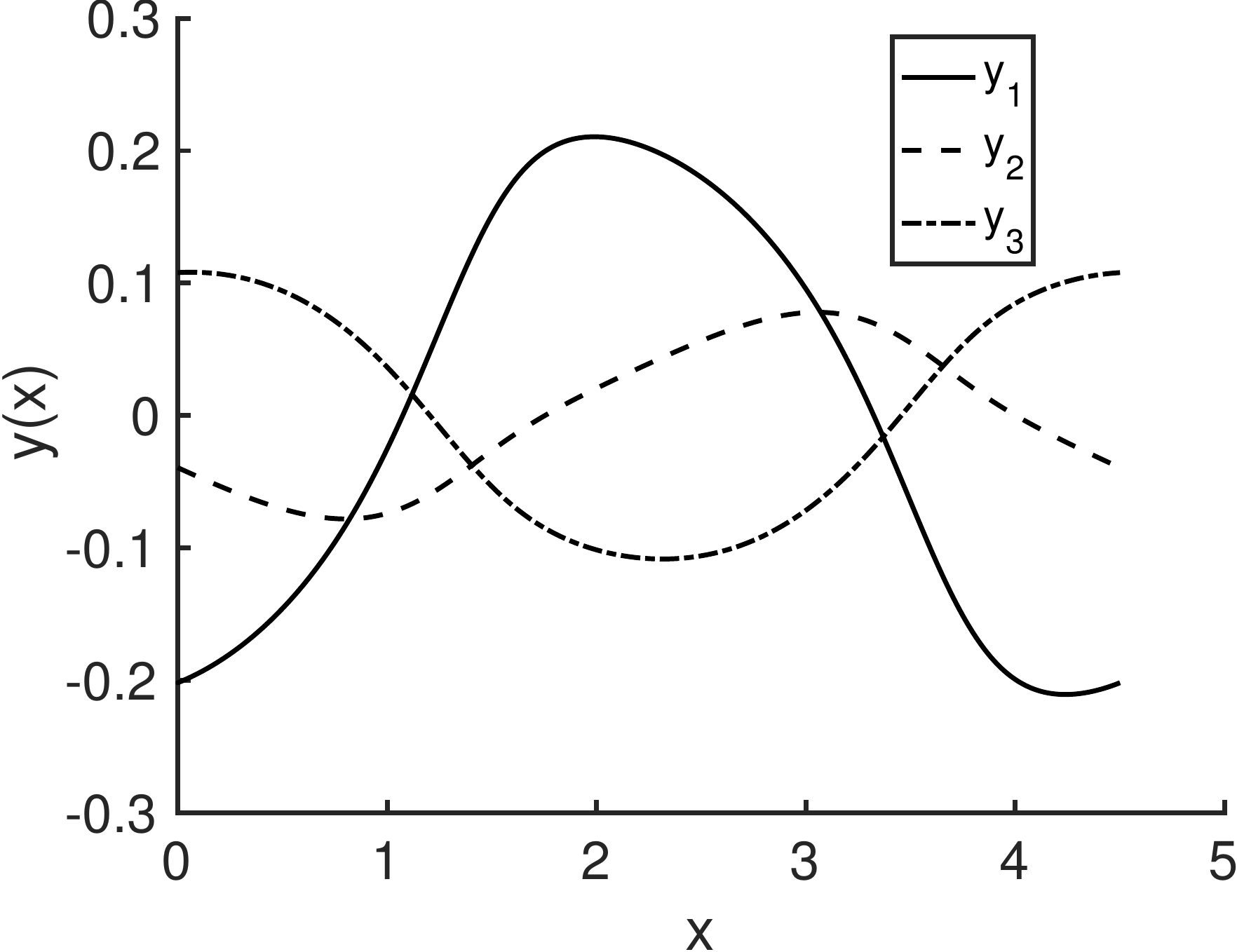} & (d) \includegraphics[scale=0.3]{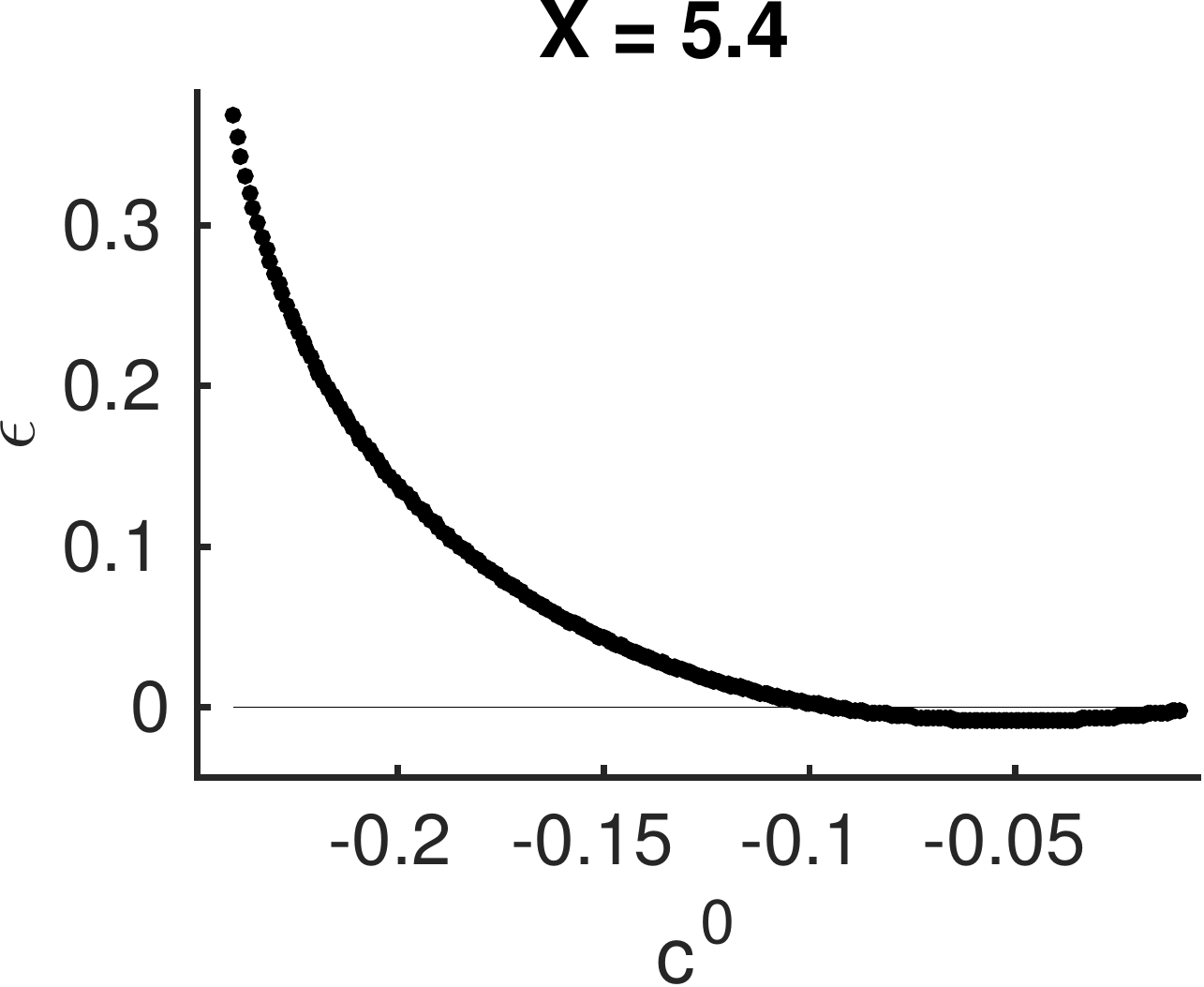}
		\end{array}
		$
	\end{center}
	\caption{(a) Stability diagram in the coordinates of shifted wave speed $c^0 = c-c_*$ and period $X$ for $\beta = -10$. Pink dots (light dots in grayscale) and black dots correspond respectively to stable and unstable waves. For a stable wave, we plot in (b) its spectrum and in (c) the wave itself, with $\beta = -10$, $c^0 = -0.3$, $X = 4.5$, and $\eps = -3.5e-3$. A dashed line marks the imaginary axis in (b). In (d) we plot a curve showing existence, up to numerical approximation, of periodic profiles of period $X = 5.4$ in the parameters $c^0$ and $\eps$ when $\beta = -10$ and the nonlinearity is cubic. A thin horizontal line marks the axis.
		}
	\label{fig126}
\end{figure}

By changing the sign of $\beta$, we find the stable periodic solutions through a sub-critical Hopf bifurcation as demonstrated in Figure \ref{fig126}. Since $\eps <0$ for the constant solution to be stable, as seen in Figure \ref{super spectrum of constant solution with nonzero c}, the periodic profile shown in Figure \ref{fig126} (c) exists through a sub-critical Hopf bifurcation. Figure \ref{fig126} (b) shows the stable spectrum of the periodic profile shown in (c). Here $\beta = -10$, $c^0 = -0.3$, $X = 4.5$ , and $\eps = -3.5e-3$. In Figure \ref{fig126} (a), we plot a stability diagram in the coordinates of shifted wave speed $c^0 = c-c_*$ and period $X$. We do not find a family of stable waves bifurcating from the Turing instability.

\subsection{Numerical stability method} To determine the spectrum of the periodic profiles, we used Hill's method. The associated eigenvalue problem is given by $Lv = \lambda v$ where the linear operator $L$ takes the form $L_{j,k}= \sum_{q=1}^{m_{jk}}f_{j,k,q}(x)\frac{\partial^q}{\partial x^q}$. The coefficients $f_{j,k,q}(x)$ are $X$ periodic. As in \cite{Deconinck2007}, we use a Fourier series to represent the coefficient functions $f_{j,k,q}$, $f_{j,k,q}(x)=\sum_{j=-\infty}^{\infty} \hat \phi_{j,k,q}e^{i2\pi jx/X}$, and write the generalized eigenfunctions as $v(x) = e^{i\xi x}\sum_{j=-\infty}^{\infty} \hat v_j e^{i\pi jx/X}$, where $\xi \in (-\pi/2X,\pi/2X]$ is the Floquet exponent. Substituting these quantities into the eigenvalue problem and equating coefficients gives an infinite dimensional eigenvalue problem for each fixed $\xi$. By truncating the Fourier series at $N$ terms and using MatLabs FFT function to determine the coefficients $\hat \phi_{j,k,q}$, we arrive at a finite dimensional eigenvalue problem $L_N^{\xi}\hat v = \lambda \hat v$, which we solve with MATLAB's eigenvalue solver. All computations were done using STABLAB \cite{STABLAB}. For further information about Hill's method and its convergence properties, see \cite{Curtis2010,Deconinck2006,Johnson2012}.

\subsection{Computational statistics}

All computations were carried out on a Macbook pro quad core or a Leopard WS desktop with 10 cores. Computing a profile took approximately 2 seconds or less, and computing the spectrum via Hill's method took on average 20-60 seconds depending on the number of modes used. We typically used 101 Floquet parameters and 41 or 81 Fourier modes when using Hill's method.  Each stability diagram took less then 24 hours to compute on the Leopard WS desktop.

\section{Discussion and open problems}\label{s:conclusions}
We have identified an analog of Turing instability occurring for $n\times n$ systems of conservation laws of
dimension $n\geq 3$, leading to a large family of spatially periodic traveling waves.
Our numerical stability investigations give convincing numerical evidence that at least some of these waves are
stable, answering the question posed in \cite{Oh2003a,Pogan2013} whether there can exist stable periodic solutions
of conservation laws.

Moreover, the same numerical investigations indicate that at least for some model parameters, the bifurcation diagram near Turing instability/Hopf bifurcation includes an open region of instability.
This opens the possibility for rigorous proof of existence of stable periodic waves through a small-amplitude bifurcation 
analysis as carried out in \cite{Mielke, Mielke1997, Sc2,SZJV} for the reaction diffusion case.
Such an analysis we consider an extremely interesting open problem.
Note, however, that it is inherently more complicated than the reaction diffusion version, involving $n+2$
bifurcation parameters $(X,c,q)$, $X,c\in \R^1$, $q\in \R^n$ rather than the two parameters of the reaction diffusion case.
For an example of intermediate complexity, we point to the recent analyses 
\cite{Matthews2000,Sukhtayev2016} of reaction diffusion equations with a single conserved quantity, featuring a
three-parameter bifurcation.


\bibliography{../auxiliary/refs}{}
\bibliographystyle{alpha}

\end{document}